\pdfoutput=1 
\documentclass{amsart}
\usepackage{xcolor}
\usepackage{hyperxmp}
\usepackage[colorlinks=false,%
			     hyperfootnotes=true,%
			     pdftex=true, %
			     implicit=true,%
			     final=true,%
			     letterpaper=true, %
			     bookmarks=true,                       
		                bookmarksnumbered=true,               
		                bookmarksopen=false,%
		pdfpagemode=UseOutlines,                
			              pdfstartview=Fit,                
			              breaklinks=true,                 
			              pageanchor=true,                 
			              plainpages=false,                
			              pdfpagelabels=true,                   
				      backref=page,
			              hyperindex=true,                      
								pdfdisplaydoctitle=true,%
								pdftitle={Computably Enumerable Equivalence Relations},%
								pdfauthor={Su Gao, Peter Gerdes},%
								pdfkeywords={computably enumerable sets,c.e.,computably enumerable equivalence relations},%
								pdfsubject={Mrclass 03D25 03D30}]{hyperref}

\newtheorem{theorem}{Theorem}[section]
\newtheorem{corollary}[theorem]{Corollary}
\newtheorem{lemma}[theorem]{Lemma}
\newtheorem{proposition}[theorem]{Proposition}

\theoremstyle{definition}
\newtheorem{definition}[theorem]{Definition}
\theoremstyle{problem}
\newtheorem{problem}[theorem]{Problem}

\usepackage{amssymb}
\newcommand{\NN}{{\mathbb N}}
\newcommand{\ZZ}{{\mathbb Z}}
\newcommand{\id}{{\rm id}}
\newcommand{\dom}{{\rm dom\,}}

\newcommand{\spec}{{\rm spec}}

\begin{document}

\title{Computably Enumerable Equivalence Relations}

\author{Su Gao}
\address{Department of Mathematics, California Institute of Technology, Pasadena,
CA 91125}
\email{{\href{mailto:sugao@its.caltech.edu}{sugao@its.caltech.edu}}}

\author{Peter Gerdes}
\address{Department of Mathematics, California Institute of Technology, Pasadena,
CA 91125}
\email{{\href{mailto:gerdes@invariant.org}{gerdes@invariant.org}}}
\thanks{The second author was supported in part by Caltech SURF 1999.}


\begin{abstract}
We study computably enumerable equivalence relations ({\em ceer}s) on $\NN$ and 
unravel a rich structural theory for a strong notion of reducibility among ceers. 
\end{abstract}

\maketitle
\tableofcontents 

\section{Introduction}

\begin{definition}
An equivalence relation $R$ on $\NN$ is a {\em ceer} if $R$ is a c.e. subset of
$\NN^2$. Such equivalence relations are also called {\em positive}.
\end{definition}

\begin{definition}
For two ceers $R_1$ and $R_2$, we say that $R_1$ is {\em m-reducible to}(or {\em
many-one reducible to}) $R_2$,
denoted $R_1\leq_mR_2$, if there is a computable function $f:\NN\to \NN$
such that, for any $x,y\in\NN$,
$$ xR_1y\Leftrightarrow f(x)R_2f(y).$$
We say that $R_1$ is {\em 1-reducible to} $R_2$, denoted $R_1\leq_1R_2$, if
there is a one-one computable function as above. Sometimes we omit the prefix and
say that $R_1$ is reducible to $R_2$, and denote $R_1\leq R_2$. When we do this, 
the reducibility is meant to be many-one. We write $R_1<_m R_2$ or $R_1<R_2$ for
the statement that $R_1\leq R_2$ but $R_2\not\leq R_1$. We write $R_1\equiv R_2$ 
if $R_1\leq R_2\leq R_1$. 
\end{definition}

Note that this reducibility is a stronger notion than the ordinary reducibility among sets.
Thus the following notion of ``completeness" is stronger than the corresponding concept
for sets.

\begin{definition}
A ceer $R$ is {\em universal} if $S\leq R$ for any ceer $S$.
\end{definition}

Positive equivalence relations were first studied by Ershov \cite{Er}, although he
did not study the reducibility notion defined above. Ershov
introduced the following notion of precompleteness.

\begin{definition}
A ceer $R$ is {\em precomplete} if $R$ has infinitely many equivalence classes
and for every partial computable function $\psi$ there is a total computable function
$f$ such that for any $n\in\NN$ with $\psi(n)\!\downarrow\,$, $f(n)R\psi(n)$.
\end{definition}

Visser \cite{Vi} studied examples of ceers in logic, especially precomplete ones.
However, it seems that Bernardi and Sorbi \cite{BS} were the first authors to
have isolated the strong reducibility notion and studied universal ceers in
the sense defined above. In particular, they proved that precomplete ceers must be
universal. Lachlan \cite{La} later showed that precomplete ceers are computably
isomorphic to each other. He also considered another natural notion of
completeness (called {\em e-complete}) and demonstrated nice properties and
natural examples of e-complete ceers.

Ceers were also studied from a different perspective by Nies \cite{Ni}. Some classes
of ceers Nies considered will be redefined in this paper and their structure of
reducibility will be investigated.

Examples of ceers in algebra, topology, logic and other areas of mathematics have
been investigated under different disguises. For example, all word problems (not
the identity problems) for semigroups and groups are in fact ceers. It is easily
derivable from the proofs of the classical unsolvability results that there are
semigroups and groups whose word problems are universal ceers. A non-trivial
example was considered by Miller \cite{Mi}, who
established that the isomorphism of finitely generated groups is a universal ceer.

Our primary goal in this paper is to develop a comprehensive theory for the structure of 
reducibility among ceers. Given the interesting circumstances that universal ceers are
in some sense better understood, we will mostly focus on non-universal ceers in this
paper. Our hope is to understand enough about non-universal ceers and the ways more
complicated ceers can be built from simpler ones, so as to eventually clarify
the structure of the strong reducibility.

The research work presented in this paper was to a large extent motivated 
by the theory of Borel reducibility
among Borel (or other definable) equivalence relations, developed by Becker, Kechris
(e.g. \cite{BK}), Friedman, Stanley (\cite{FS}) and others. Although there are no
direct connections between the two subjects, the general methodology we employ
in the current research is much similar to the theory of Borel equivalence relations.

Let us briefly mention some speculations about the advantages to consider the
strong reducibility notion. First of all, the most interesting applications of
the theory of equivalence relations are often classification problems for various
kinds of structures in mathematics. Furthermore, to classify a structure is usually
meant to obtain a complete set of invariants for the structure. What we are investigating
about ceers under the strong reducibility can then be viewed as an abstract framework
for a study of the possibilities to effectively {\em compute} complete invariants. In this
sense our reducibility is arguably more natural than the reducibility among sets.

Second, a ceer can be viewed as a c.e. decomposition of the universe. Hence the 
reducibility among ceers as equivalence relations generalizes 
various notion of simultaneous reduction for sequences of sets or reduction 
for ``promise" problems. Also, natural non-universal ceers are easy to construct.
This is in contrast with the situation in computability theory for sets, where
sophisticated constructions are usually needed to ensure existence. 

The rest of the paper is organized as follows. Section 2 contains some preliminaries
about ceers. In sections 3-6 we consider ceers with simplest formations and explore
the reducibility among them. In sections 7 and 8 we defined two notion of jump operations
for ceers. As we have mentioned above, more complicated ceers can thus be formed and
more interesting properties can be obtained. In sections 8 and 9 we clarify the
structure of reducibility among the ceers constructed in earlier sections. Section 10
contains a summary of the theory and some open problems to consider in future researches.

\section{Indices for ceers}
Our basic notation about partial computable functions and c.e. sets follows
that of \cite{Ro}.
The partial computable function with index $e$ is denoted by $\varphi_e$.
The c.e. set with index $e$ is denoted by $W_e$.
If $\psi$ is a partial computable function and $x\in\NN$, then the domain of 
$\psi$ is denoted by $\dom\psi$; $\psi(x)\!\downarrow\,$
indicates that $x\in\dom\psi$ and $\psi(x)\!\uparrow\,$ otherwise.
For a partial computable function $\psi$ and $n\in\NN$, $\psi^n$ denotes the
partial computable function obtained by $n$ iterated compositions of $\psi$.
This can be defined precisely by induction on $n$: $\psi^0(x)=x$ 
and $\psi^{n+1}(x)=\psi(\psi^n(x))$ for any $x\in\NN$.

Throughout the paper we adopt the following notation. Denote
$K=\{ x\in\NN\,|\, \varphi_x(x)\!\downarrow\,\}$ and $\overline{K}=\NN\setminus K$. 
For $i\in\NN$, let $K_i=\{ x\in K\,|\,\varphi_x(x)\!\downarrow\,=i\}$. We use the
notation $\langle \cdot,\cdot\rangle$ to denote a fixed computable bijection
between $\NN^2$ and $\NN$ and call it the {\em coding function}. Coding functions
of arbitrary arities are similarly represented.

Next we fix some notation and note some basic facts about ceers.

For a ceer $R$ on $\NN$ and $x\in\NN$, denote by $[x]_R$ the $R$-equivalence class
containing $x$ and by $|[x]_R|$ the cardinality of $[x]_R$.

The set of all ceers is computably enumerable. Let $R_e$ be the equivalence 
relation generated by the set $\{(x,y)\,|\,\langle x,y\rangle
\in W_e\}$. Then $R_e$ is a ceer and all ceers appear in this enumeration as some $R_e$. 
We will regard this as the canonical enumeration for all ceers and say that 
$e$ is a {\em canonical c.e. index} for the relation $R_e$. 

Let $R_\infty$ be defined by
$$ \langle x,z\rangle R_\infty \langle y,z\rangle\Leftrightarrow xR_z y$$
for $x,y,z\in\NN$. (It is to be understood in this definition, as well as in
all definitions of a similar form in the sequel, that no other pair of distinct
numbers is in the relation.) It is obvious that $R_\infty$ is a
universal ceer.

There are other procedures which can generate all ceers. The following
definition is due to Ershov.

\begin{definition} For each partial computable function $f$, define the ceer
$\eta_f$ by
$$ x\eta_fy\Leftrightarrow \exists n,m (f^n(x)\!\downarrow\, =f^m(y)\!\downarrow\, ),
\mbox{ for $x,y\in\NN$}.$$
Then $\eta_f$ is the ceer
generated by the graph of $f$.
If $e$ is an index for $f$, then we also write $\eta_f$ as $\eta_e$, and we
call $e$ the {\em iterative c.e. index} for $\eta_e$.
\end{definition}

Ershov \cite{Er} showed that every ceer has an iterative c.e. index. We note
below that the two indexing systems are essentially the same.

\begin{theorem}
There is a computable isomorphism $\rho$ such that $\eta_e=R_{\rho(e)}$.
\end{theorem}
\begin{proof}
A proof of Ershov's result indicates, in fact it is not hard to see directly, 
that there is a one-one
computable function $f(n,e)$ such that $R_e=\eta_{f(n,e)}$ for any $n,e\in\NN$.
Similarly, there is a one-one computable function $g(n,e)$ such that
$\eta_e=R_{g(n,e)}$ for any $n,e\in\NN$. With this padding property it is then
routine to define a computable isomorphism
$\rho$ by a back-and-forth construction as follows.

Let $\rho(0)=f(0,0)$. If $f(0,0)=0$ then $\rho^{-1}(0)$ is already defined.
Otherwise, we define $\rho^{-1}(0)$ to be $g(k,0)$ where $k$ is the least
such that $g(k,0)\neq 0$. In general suppose both $\rho$ and $\rho^{-1}$ are
defined for $i<n$. If $n\in\{ \rho^{-1}(0),\dots,\rho^{-1}(n-1)\}$ then
$\rho(n)$ is already defined. Otherwise let $\rho(n)$ be $f(k,n)$ where $k$
is the least such that $f(k,n)\not\in\{ \rho(0),\dots,\rho(n-1)\}$. Then 
in a similar way define $\rho^{-1}(n)$. The resulting $\rho$ is
one-one, onto and computable. 
\end{proof}

Ceers can also be generated by computable actions of countable computable groups
on $\NN$. When the ceer has no finite equivalence classes, the group can be
taken to be the additive group of $\ZZ$.

It will be useful to note that the reducibility among ceers is a $\Sigma_3$ relation.
Specifically, the index set $\{ (e,i)\,|\, R_e\leq R_i\}$ is $\Sigma_3$. This is by
a straightforward computation. It follows that the set $\{ (e,i)\,|\, R_e\equiv R_i\}$
is also $\Sigma_3$. In the subsequent sections we will show that both 
reducibility and bi-reducibility among ceers are $\Sigma_3$-complete.

In the subsequent sections we will consider various classes of ceers, and we will
adopt the following terminology.

\begin{definition} Let $\mathcal C$ be a non-empty class of ceers. A ceer $R$ is
said to be {\em essentially} in $\mathcal C$ if $R\leq S$ for some $S\in {\mathcal C}$.
$R$ is {\em universal} for $\mathcal C$ if $R\in{\mathcal C}$ and $S\leq R$ for any
$S\in{\mathcal C}$. If $\mathcal C$ has a universal relation in it, then $R$ is
{\em essentially universal} for $\mathcal C$ if $R\equiv S$ for some $S$ universal
for $\mathcal C$.
\end{definition}

Much of our investigation is focused on finding and describing universal relations for
different classes of ceers.

\section{Computable equivalence relations and partial classifiability}
Computable equivalence relations are the simplest ones among ceers with respect
to reducibility. 

\begin{definition}
For each $n\in \NN$ and $n>0$, let $\id(n)$ denote the following equivalence
relation on $\NN$ with $n$ equivalence classes:
$$ x\id(n)y\Leftrightarrow {x\equiv y} (\mbox{mod}\,n), \mbox{ for $x,y\in \NN$}.$$
Let $\id(\NN)$ or $\omega$ be the identity relation on $\NN$:
$$ x\omega y\Leftrightarrow x=y, \mbox{ for $x,y\in\NN$}.$$
\end{definition}

\begin{definition}
Let $R$ be a ceer on $\NN$. A subset $T\subseteq \NN$ is called a {\em transversal}
for $R$ if $T$ meets each $R$-equivalence class in exactly one element.
\end{definition}

Then we have the following complete characterization and classification of computable
equivalence relations.

\begin{proposition} $~$ \
\begin{enumerate}
\item[(1)] For each $n\in \NN$ and $n>0$, $\id(n)$ is computable; and $\omega$ is computable.
\item[(2)] $\id(1)<\id(2)<\cdots<\id(n)<\cdots<\omega$.
\item[(3)] If $R<\omega$, then there is some $n\in\NN$ such that $R\equiv \id(n)$.
\end{enumerate}
\end{proposition}

\begin{proof} (1) and (2) are obvious. For (3), suppose $R<\omega$, then it follows that
$R$ is computable. If $R$ has infinitely many equivalence classes, then $\omega\leq R$
via the following computable function $f$ defined by primitive recursion:
$$\begin{array}{rcl}
f(0)& = & 0 \\
f(x+1)& = & \mu y(y>f(x)\land \forall z\leq x\,\lnot(yRf(z))).
\end{array}
$$
It follows that $R$ has only finitely many equivalence classes. Let $n\in\NN$ be the number
of equivalence classes of $R$. Then it is easy to see that $R\equiv \id(n)$.
\end{proof}

\begin{proposition}
The following are equivalent for any ceer $R$:
\begin{enumerate}
\item[($i$)] $R$ is computable;
\item[($ii$)] $R\leq\omega$;
\item[($iii$)] $R$ has an c.e. transversal;
\item[($iv$)] $R$ has a computable transversal;
\item[($v$)] $T(R)=_{\rm def}\{ \min(C) \,|\, \mbox{$C$ is an $R$-equivalence class}\}$ 
is computable.
\end{enumerate}
\end{proposition}
\begin{proof}
Note that 
$(i)\Rightarrow (v)\Rightarrow (iv)\Rightarrow (iii)$
and $(ii)\Rightarrow (i)$ are obvious. We are only left to prove $(iii)\Rightarrow (ii)$. 
Suppose $T$ is an c.e. transversal for $R$.
Consider the function $g:\NN\to T$ defined by letting $g(x)$ be the unique element
$y\in T$ with $xRy$. Then $g$ is total because $T$ is a transversal. It is obvious that
the graph of $g$ is c.e., hence $g$ is computable. But then $g$ is a reduction from $R$ to $\omega$.
\end{proof}

The proposition implies that computable ceers are exactly the kind of ceers for which
natural numbers can be effectively assigned to equivalence classes as complete 
invariants. Such equivalence relations should be considered classifiable in an ideal
sense. Historically people have also considered the possibility of assigning
an infinite system of natural numbers as complete invariants (e.g. \cite{Ra}),
in the following we remark that infinite systems of invariants do not help.

\begin{proposition}
Let $R$ be a non-computable ceer. Then there is no computable function $f(m,n)$
such that, for any $x,y\in\NN$,
$$ xRy\Leftrightarrow \forall m (f(m,x)=f(m,y)). $$
\end{proposition}
\begin{proof}
The existence of such an $f$ would imply that $R$ is co-c.e.; since $R$ is c.e.
it would follow that $R$ is computable, contradiction.
\end{proof}

Next we formulate a notion of partial classifiability for ceers.

\begin{definition} For each partial computable function $\psi$ define ceer $P_\psi$ by
$$ xP_\psi y\Leftrightarrow x=y \lor \psi(x)\!\downarrow\, =\psi(y)\!\downarrow\,, \mbox{ for $x,y\in\NN$.}$$
Since $\psi$ can be viewed as a function partially assigning
complete invariants for (at least) the non-trivial equivalence classes, we call ceers
of the form $P_\psi$ {\em partially classifiable} or in abbreviated form, {\em PC}. 
\end{definition}

It turns out that there is a universal PC relation.

\begin{definition}
The ceer $H$ (standing for {\em halting equivalence}) is defined by
$$ xHy\Leftrightarrow x=y\lor \varphi_x(x)\!\downarrow\, =\varphi_y(y)\!\downarrow\, $$
for $x,y\in\NN$.
\end{definition}

It is immediate from the definition that $H$ is PC.

\begin{proposition}
$H$ is universal for PC relations, i.e., for any PC relation $R$, $R\leq H$.
\end{proposition}
\begin{proof}
Suppose $R$ is PC, as witnessed by partial computable function $\psi$.
Let $f$ be one-one computable function obtained by $s$-$m$-$n$ Theorem such that
$\varphi_{f(x)}(y)=\psi(x)$. Then $f$ witnesses that $R\leq H$.
\end{proof}

The following simple observation will be useful in {\S} 7.
\begin{proposition}
Let $E$ be a ceer without computable classes. Let $R$ be a PC-relation. Then
$E\not\leq R$.
\end{proposition}
\begin{proof}
Suppose $f$ is a computable function witnessing $E\leq R$. Suppose $\psi$ is
a partial computable function with $R=P_\psi$. Then we claim that, for every
$x\in \NN$, $\psi(f(x))$ is defined. Otherwise, we have that $[f(x)]_R$ is
a singleton containing only $f(x)$, thus by the assumption that $f$ is a reduction,
we would have that $[x]_E$ is computable. Therefore we in fact have that
$\psi\circ f$ is a total computable function; moreover, $\psi\circ f$ witnesses
that $E\leq\omega$, again implying that every $E$-class is computable, contradiction.
\end{proof}

\section{Finite dimensional ceers}

\begin{definition} Let $n\in\NN$ and $n>0$.
A ceer $R$ is {\em $n$-dimensional} if there are pairwise disjoint c.e. sets
$A_1,\dots,A_n$ such that
$$ xRy\Leftrightarrow x=y \lor \exists i\leq n (x,y\in A_i).$$
We denote the above relation by $R_{A_1,\dots,A_n}$.
A ceer is {\em finite dimensional} if it is $n$-dimensional for some $n\in\NN$.
\end{definition}

We have the following easy facts about $1$-dimensional ceers.
\begin{proposition} Let $A, B$ be non-empty c.e. sets. 
\begin{enumerate}
\item[(1)] $R_A$ is computable iff $A$ is computable.
\item[(2)] $A\leq_1B$ iff $R_A\leq_1 R_B$.
\item[(3)] If $R_A\leq R_B$, then $A\leq_mB$.
\end{enumerate}
\end{proposition}
\begin{proof}
For (1) one direction follows immediately from the definition.
Conversely, suppose $R_A$ is computable, then it has a computable transversal $T$. Note that $T\cap A$ is a singleton, therefore $T\setminus A$ is computable. But $T\setminus A$ is just the complement of $A$, it follows that $A$ is computable. (2) and (3) are obvious.
\end{proof}

The proposition indicates that the reducibility among $1$-dimensional relations is very 
similar to that among c.e. sets. In fact, the proof of (1) shows that $R_A$ (as a subset of 
$\NN^2$) is in the same Turing degree as $A$. Because of the transparency of the reducibility
among $1$-dimensional ceers modulo that among the c.e. sets, we consider the $1$-dimensional
ceers completely classified.

\begin{proposition}
Let $n\in\NN$ and $n>0$.
\begin{enumerate}
\item[(1)] If $R$ is an $n$-dimensional ceer, then $R\leq_1 R_{K_0,\dots,K_{n-1}}$.
\item[(2)] $R_{K_0,\dots,K_{n-1}}<R_{K_0,\dots,K_n}$.
\item[(3)] Let $k\leq n$ and $R=R_{A_1,\dots,A_n}$.
Then $R$ is essentially $k$-dimensional iff at most $k$ of the sets $A_1,
\dots,A_n$ are non-computable.
\end{enumerate}
\end{proposition}
\begin{proof}
First consider (1). Let $A_0,\dots,A_{n-1}$ be disjoint c.e. sets. We are to define a reduction
to witness $R_{A_0,\dots,A_{n-1}}\leq R_{K_0,\dots,K_{n-1}}$. First define a
partial computable function $F$ as follows.
$$ F(e,x,y)=\left\{\begin{array}{ll}
i & \mbox{ if $i<n, x\in A_i$ and $y=\varphi_e(x)$} \\
\!\uparrow\, & \mbox{ otherwise}
\end{array}\right.
$$
Let $s(e,x)$ be a one-one total computable function with
$F(e,x,y)=\varphi_{s(e,x)}(y)$. Then by the fixed point theorem, let $e_0$
be such that $\varphi_{e_0}(x)=s(e_0,x)$, for all $x$. Finally let $f(x)=\varphi_{e_0}(x)$.
Then $f$ is one-one computable. We verify that $f$ is the required reduction.
If $x\in A_i$ for some $i<n$, then $F(e_0,x,f(x))=i$ by definition of $F$,
which means $\varphi_{s(e_0,x)}(f(x))=\varphi_{f(x)}(f(x))=i$, or $f(x)\in K_i$.
If $x\not\in A_i$ for any $i<n$, then by the definition of $F$, $F(e_0,x,f(x))\!\uparrow\,$,
which means $f(x)\not\in K_i$ for any $i<n$.

For (2), note first that it follows from (1) that 
$R_{K_0,\dots,K_{n-1}}\leq R_{K_0,\dots,K_n}$. To see that
$R_{K_0,\dots,K_n}\not\leq R_{K_0,\dots,K_{n-1}}$, just note that for any reduction
function $h$ there is some $i\leq n$ such that $h$ sends the class $K_i$ to a single
element. This implies that $K_i$ is computable, a contradiction.

(3) can be shown by a similar argument as in the last paragraph.
\end{proof}

The proposition establishes a hierarchy for finite dimensional ceers with a universal
relation in each level.

At this point it is natural to ask whether the identity relation is reducible to
all non-computable ceers. We prove that this is not the case. Furthermore, we
completely characterize those finite dimensional ceers to which the identity
relation is reducible. 
Recall the definition of simple sets by Post.
\begin{definition} An c.e. set $A\subseteq \NN$ is {\em simple} if 
\begin{enumerate}
\item[(a)] the complement of $A$ is infinite, and
\item[(b)] for any infinite c.e. set $B\subseteq\NN$, $A\cap B\neq\emptyset$.
\end{enumerate}
\end{definition}

It is well known that simple sets exist and that none of them is $m$-complete among the c.e. sets.
\begin{proposition} Let $A_1,\dots,A_n$ be disjoint c.e. sets the complement of whose
union is infinite. Then
$$ \omega\leq R_{A_1,\dots,A_n}\Leftrightarrow \mbox{ $A_1\cup\dots\cup A_n$ is not simple}.$$
\end{proposition}
\begin{proof} Let $A=\bigcup_{i\leq n}A_i$.
Suppose $f$ witness the reduction $\omega\leq R_{A_1,\dots,A_n}$. Let $S$ be the range of $f$. Then $S\cap A$ contains at most $n$ elements. It follows that $S\setminus A$ is an infinite c.e. set. Therefore $A$ is not simple. Conversely, if $A$ is not simple, then its complement contains an infinite c.e. set. Let $f$ be a computable function enumerating this set. Then $f$ is a reduction from $\omega$ to $R_{A_1,\dots,A_n}$. 
\end{proof}

Next let us remark that all finite dimensional ceers are PC.

\begin{proposition} All finite dimensional ceers are PC.
\end{proposition}
\begin{proof}
Let $A_1,\dots,A_n$ be disjoint c.e. sets. Then $R_{A_1,\dots,A_n}$ is PC as witnessed
by the following partial computable function $\psi$, defined by
$$ \psi(x)=\left\{\begin{array}{ll}
i & \mbox{ if $x\in A_i$} \\
\!\uparrow\, & \mbox{ otherwise}
\end{array}\right.
$$
\end{proof}

Thus the halting equivalence $H$ is an upper bound for all finite dimensional
ceers in the reducibility order. In fact, the definitions suggest that $H$ is
quite a natural limit for all $R_{K_0,\cdots,K_n}$. Nevertheless, it is not clear
whether $H$ is a {\em least} upper bound.

We are now ready to show that reducibility and bi-reducibility
among ceers are $\Sigma_3$-complete.

\begin{proposition} The set $\{ e\,|\, R_e\leq \omega\}$ is $\Sigma_3$-complete.
\end{proposition}
\begin{proof}
It is well known that the set Comp$=\{ i\,|\, W_i\mbox{ is computable}\}$ is $\Sigma_3$-complete
(see \cite{Ro} Theorem 14-XVI). For each $i\in\NN$, consider the one-dimensional ceer $R_{W_i}$.
It is obvious that $R_{W_i}\leq\omega$ iff $R_{W_i}$ is computable iff $W_i$ is computable.
Moreover, it is easy to see that there is a computable function $f$ such that 
$R_{f(i)}=R_{W_i}$. Therefore Comp is reducible to the set in question.
\end{proof}

This immediately gives the $\Sigma_3$-completeness of reducibility and bi-reducibility.

\begin{corollary}
The reducibility and bi-reducibility among ceers is $\Sigma_3$-complete.
\end{corollary}
\begin{proof}
The assertion about reducibility follows trivially from the above proposition.
For bi-reducibility consider, for each $i\in\NN$, the one-dimensional ceer
$R_{2W_i}$, where $2W_i=\{ 2x\,|\, x\in W_i\}$.
Then $\omega\leq R_{2W_i}$ for every $i\in\NN$. Moreover, $R_{2W_i}\leq \omega$ iff
$2W_i$ is computable iff $W_i$ is computable. We thus
have that the set $\{ e| R_e\equiv \omega\}$ is $\Sigma_3$-complete.
\end{proof}

Next we consider the complexity of various other sets arising from the
reducibility of ceers.

\begin{lemma}
The index set $\{ e\,|\, \exists i (\mbox{ $W_i$ is infinite and } W_i\cap W_e=\emptyset )\}$
is $\Sigma_3$-complete.
\end{lemma}
\begin{proof}
It is well known that the set $\{ e\,|\, W_e \mbox{ is simple}\}$ is $\Pi_3$-complete.
A proof of this can be found in \cite{Ro} Exercise 14-31, which uses eventually a theorem
of Dekker (\cite{Ro} Theorem 9-XVI). Note that Dekker's proof effectively shows
that there is a computable function $f$ such that, for any $x$, $W_{f(x)}$ is co-infinite,
and $W_x$ is not computable iff $W_{f(x)}$ is simple (in fact hypersimple). This reduction
function $f$ witnesses that Comp is reducible to our set in question.
\end{proof}

\begin{proposition}
The set $\{ e\,|\, \omega\leq R_e\}$ is $\Sigma_3$-complete.
\end{proposition}
\begin{proof}
This now follows easily from Proposition 4.5 and Lemma 4.9.
\end{proof}

By similar proofs we obtain the following results.
\begin{proposition}
The following index sets are $\Sigma_3$-complete:
\begin{enumerate}
\item[(1)] $\{ e\,|\, R_e\leq R_A\}$ and $\{ e\,|\, R_e\equiv R_A\}$,
where $A\subseteq \NN$ is non-empty, c.e., not simple and has infinite complement.
\item[(2)] $\{ e\,|\, R_e\mbox{ is essentially $n$-dimensional}\}$, where $n>0$.
\item[(3)] $\{ e\,|\, R_e\mbox{ is essentially universal $n$-dimensional}\}$, 
where $n>0$.
\item[(4)] $\{ e\,|\, R_e\mbox{ is essentially finite dimensional}\}$.
\end{enumerate}
\end{proposition}
\begin{proof} Let us sketch the proof for (1); the other parts are similar.
Let $A$ be as given and $C$ be a computable infinite subset of $\NN\setminus A$.
Fix a computable bijection $\theta$ from $\NN$ onto $C$. For each $e$ consider the
$2$-dimensional ceer $R_{A,\theta(W_e)}$. By Proposition 4.3 it is easy to see
that $R_{A,\theta(W_e)}\leq R_A$ iff $W_e$ is computable, and when this happens
we also have $R_{A,\theta(W_e)}\equiv R_A$. Similar to the proof of Proposition 4.7 
we can then construct reductions from Comp to the index sets 
$\{ e\,|\, R_e\leq R_A\}$ and $\{ e\,|\, R_e\equiv R_A\}$.
\end{proof}

\section{FC relations and CC relations}

Since ceers with only finitely many equivalence classes are not very interesting in terms of reducibility, we call such relations {\em trivial}.
\begin{definition}
A ceer $R$ is {\em FC} (standing for {\em finite classes}) if all of its equivalence 
classes are finite.
$R$ is {\em CC} (for {\em computable classes}) if all of its equivalence classes are 
computable.
\end{definition}

Let us summarize some basic facts.
\begin{proposition} $~$\
\begin{enumerate}
\item[(1)] All FC relations are CC.
\item[(2)] If $R_1\leq R_2$ and $R_2$ is CC, then $R_1$ is CC.
\item[(3)] If $R$ is a PC-relation whose every equivalence class is non-trivial, then
$R$ is computable. Therefore, there are CC relations which are not PC. In fact, for each
non-computable PC-relation $R$ there is $S\equiv R$ such that $S$ is not PC.
\item[(4)] For each FC relation $R$ there is $S\equiv R$ such that $S$ is not FC.
\end{enumerate}
\end{proposition}
\begin{proof}
(1),(2) and the first part of (3) are obvious. For the second part of (3) and (4) define
$S$ by
$$ \langle x,u\rangle S\langle y,v\rangle \Leftrightarrow xRy, $$
then every equivalence class of $S$ is infinite, and it is easy to see that $S\equiv R$.
\end{proof}

We have the following non-reducibility result between finite dimensional relations and 
CC relations.
\begin{proposition}
Let $R$ be a non-computable finite dimensional ceer and $S$ be a non-computable CC relation.
Then $R\not\leq S$ and $S\not\leq R$.
\end{proposition}

\begin{proof}
Note the general fact: if $R_1\leq R_2$ via $f$ and $[y]_{R_2}$ is computable, then
the preimage set $C=\{ x\,|\, f(x)R_2y\}$ is also computable. In fact, either $C$ is empty
or $C$ is a single $R_1$-equivalence class. Moreover $C\leq_m [y]_{R_2}$ via $f$. 
Let $R$ and $S$ be given as in the hypothesis. Then
$R\not\leq S$ is an immediate consequence of the above fact. 
For $S\not\leq R$, suppose $R=R_{A_1,\dots,A_n}$ where $A_1,\dots,A_n$ are disjoint
c.e. non-computable. Assume that $f$ witnesses $S\leq R$ and suppose $B=\{ f(x)\,|\, x\in\NN\}$.
Let $A=\bigcup_{i\leq n}A_i$. If $B\cap A=\emptyset$, then $f$ is in fact a witness for $S\leq\omega$, contradicting the hypothesis that $S$ is non-computable.
Otherwise, let $C_1,\dots,C_k$ ($k\leq n$) be the $S$-equivalence classes with 
$f(x)\in A_i$ for $i\leq k$ and $x\in C_i$. Choose arbitrary elements $y_i\in A_i$ for $i\leq k$
and define
$$ g(x)=\left\{\begin{array}{ll}
y_i & \mbox{ if $i\leq k$ and $x\in C_i$} \\
f(x) & \mbox{ otherwise}
\end{array}\right.
$$
Then $g$ witnesses that $S\leq \omega$, again contradicting the hypothesis that $S$ is
non-computable.
\end{proof}

We next consider the question of when $\omega\leq R$ for an FC relation $R$. 
\begin{definition}
For an c.e. set $A\subseteq\NN$, let $F_A$ be the equivalence relation
defined by
$$ xF_Ay\Leftrightarrow x=y\lor \forall z (\min(x,y)\leq z\leq \max(x,y)
\to z\in A), \mbox{ for $x,y\in\NN$}.$$
Then $F_A$ is a ceer. 
\end{definition}

If $A$ is an c.e. set whose complement is infinite, then $F_A$ is an FC relation.
In analogy with the result of preceding section we have the following proposition. Let us
recall the definition of hypersimple sets, also due to Post.

\begin{definition} $~$\
\begin{enumerate}
\item[(a)] Let $A\subseteq \NN$ be an infinite set and $f$ a total function on $\NN$. We say
that $f$ {\em majorizes} $A$ if, letting $z_0,z_1,\dots$ be the enumeration of $A$ in strictly
increasing order, we have that $(\forall n) f(n)\geq z_n$.
\item[(b)] An c.e. set $A\subseteq\NN$ is {\em hypersimple} if the complement of $A$ is infinite
and there is no computable $f$ such that $f$ majorizes $\NN\setminus A$.
\end{enumerate}
\end{definition}

\begin{proposition}
Let $A\subseteq\NN$ be an c.e. set whose complement is infinite. Then
$$ \omega\leq F_A\Leftrightarrow \mbox{ $A$ is not hypersimple}.$$
\end{proposition}
\begin{proof}
Suppose $f$ witnesses that $\omega\leq F_A$. Then define by recursion
a function $g$ as follows.
$$\begin{array}{rcl}
g(0) & = & f(0), \\
g(n+1) & = & f(x), \mbox{ where $x$ is the least with $f(x)>g(n)$.}
\end{array}
$$
Since the complement of $A$ is infinite, $g$ is total computable. Letting $h(n)=g(n+1)$,
then $h$ majorizes $\NN\setminus A$, therefore $A$ is not hypersimple. Conversely, if
$A$ is not hypersimple, let $h$ be a computable function majorizing $\NN\setminus A$
and let $z_0,z_1,\dots$ be an enumeration of $\NN\setminus A$ in strictly increasing order.
We define a computable function $f$ enumerating an infinite set of pairwise inequivalent
elements of $F_A$, then $f$ witnesses that $\omega\leq F_A$.
Let 
$$\begin{array}{rcl}
f(0)&=&0, \\
f(n+1)&=&h(f(n)).
\end{array}
$$ 
Then $f(n)\leq z_{f(n)}<h(f(n))=f(n+1)$. Since
$z_{f(n)}\not\in A$, by definition of $F_A$ we have that $f(n)$ and $f(n+1)$ are inequivalent.
\end{proof}

The following results concern the complexity of being FC, essentially FC and CC.

\begin{lemma}
There is a computable function $f$ such that, for any $x\in\NN$,
$$ \begin{array}{rcl}
\mbox{$W_x$ is co-finite}& \Rightarrow & \mbox{$R_{f(x)}$ is not CC, and} \\
\mbox{$W_x$ is co-infinite}& \Rightarrow &\mbox{$R_{f(x)}$ is FC.}
\end{array}
$$
\end{lemma}
\begin{proof}
Given an c.e. set $W$ define a ceer $R^W$ by
$$ xRy\Leftrightarrow x=y\lor (\forall z (\min(x,y)\leq z\leq\max(x,y)\to z\in W)
\land x,y\in K). $$
Then if $W$ is co-finite $R$ is a finite modification of $R_K$, hence $R$ is not CC.
If $W$ is co-infinite then $R$ is obviously FC.
\end{proof}

\begin{proposition}$~$\
\begin{enumerate}
\item[(1)] The set $\{ e\,|\, \mbox{$R_e$ is FC}\}$ is $\Pi_3$-complete.
\item[(2)] The set $\{ e\,|\, R_e\mbox{ is essentially FC}\}$ is $\Pi_3$-hard.
\end{enumerate}
\end{proposition}
\begin{proof} These follow from the preceding lemma and some additional observations.
First of all, note that the set $\{ e\,|\, \mbox{$W_e$ is co-finite}\}$
is known to be $\Sigma_3$-complete.
Then, for (1) just note that being FC is $\Pi_3$. For (2) note that if $R$ is not CC, then
$R$ is not essentially FC. 
\end{proof}

\begin{lemma}
The set $\{ e\,|\, \forall x \mbox{ $\{ y\,|\, \langle x,y\rangle\in W_e\}$
is computable}\}$ is $\Pi_4$-complete.
\end{lemma}
\begin{proof}
Denote the set by $X$. $X$ is obviously $\Pi_4$. 
To see that it is $\Pi_4$-complete, let $A$ be an
arbitrary $\Pi_4$ set and let $B\in\Sigma_3$ be such that
$$ x\in A\Leftrightarrow \forall y (\langle x,y\rangle\in B), \mbox{ for all $x$.}$$
Let $f$ be a reduction function from $B$ to Comp. Define a computable function $g$ such
that, for any $e,x,y\in\NN$,
$$ \langle y,x\rangle\in W_{g(e)}\Leftrightarrow x\in W_{f(\langle e,y\rangle)}.$$
Then for any $e\in\NN$,
$$\begin{array}{rcl}
e\in A &\Leftrightarrow & \forall y (\langle e,y\rangle\in B) \\
&\Leftrightarrow & \forall y (f(\langle e,y\rangle)\in\mbox{Comp}) \\
&\Leftrightarrow & \forall y (W_{f(\langle e,y\rangle)} \mbox{ is computable}) \\
&\Leftrightarrow & \forall y (\mbox{ the set $\{ x\,|\, \langle y,x\rangle\in W_{g(e)}\}$
is computable}) \\
&\Leftrightarrow & g(e)\in X
\end{array}
$$
This shows that $A\leq_m X$. Since $A$ is arbitrary, $X$ is $\Pi_4$-complete.
\end{proof}

\begin{proposition}
The set $\{ e\,|\, \mbox{$R_e$ is CC}\}$ is $\Pi_4$-complete.
\end{proposition}
\begin{proof} A straightforward computation shows that it is $\Pi_4$. To see that
it is $\Pi_4$-complete, for any c.e. set $W$ define a ceer $R^W$ by
$$ \langle x,y\rangle R^W\langle x,z\rangle\Leftrightarrow y=z\lor (\langle x,y\rangle,
\langle x,z\rangle\in W) $$
for any $x,y,z\in \NN$. Let $X$ be the index set in the preceding lemma.
Then $e\in X$ iff the relation $R^{W_e}$ is CC.
\end{proof}

\begin{corollary}$~$\
\begin{enumerate}
\item[(1)] There is no universal FC relation.
\item[(2)] There is no universal CC relation.
\item[(3)] There is an CC relation which is not essentially FC.
\end{enumerate}
\end{corollary}
\begin{proof}
For (1) and (2) note that if there were a universal FC or CC relation, 
then being essentially FC or being CC would be $\Sigma_3$, contradicting
Proposition 5.8(2) and Proposition 5.10. For (3),
if every CC relation were essentially FC, then the sets 
$\{ e\,|\, \mbox{$R_e$ is CC}\}$ and
$\{ e\,|\, \mbox{$R_e$ is essentially FC}\}$ would be the same. 
But this is impossible since
the former is $\Pi_4$-complete and the latter is $\Sigma_4$.
\end{proof}

Part (3) of this corollary asserts the abstract existence of some ceers of 
which we do not yet have any constructed example.

\section{Bounded relations}

\begin{definition}
For $k\in\NN$ and $k>0$, we say that $R$ is {\em $k$-bounded} if every $R$-equivalence class contains at most $k$ elements. $R$ is {\em bounded} if it is
$k$-bounded for some $k\in\NN$. 
\end{definition}

The following theorem shows that $\omega\leq R$ for all bounded relations $R$ but there is no uniform effective way to find the reduction.

\begin{theorem} $~$\
\begin{enumerate}
\item[(1)] If $R$ is a bounded relation, then $\omega\leq R$.
\item[(2)] Let $k\in\NN$ and $k>1$. There is no partial computable function $\rho(e,x)$ such that, whenever $R_e$ is a $k$-bounded ceer then $\rho(e,x)\!\downarrow\,$, for all $x\in\NN$ and $\rho(e,\cdot)$ witnesses the reduction from $\omega$ to $R_e$.
\end{enumerate}
\end{theorem}

\begin{proof} (1)
Suppose $R$ is $k$-bounded. For each $i\leq k$, let $N_i$ be the number of $R$-equivalence classes containing exactly $i$ many elements. Since
$\sum_{i\leq k}N_i=\infty$, there is a biggest $l\leq k$ such that $N_l=\infty$.
Let $F$ be the set of all elements whose $R$-equivalence class has more than $l$ elements. Then $F$ is a finite set. Now let $\Gamma$ be an enumeration of the pairs of elements in $R$. We define a computable reduction from $\omega$ to $R$ by listing an infinite set of pairwise non-$R$-equivalent elements by stages. At stage $n$ suppose a finite set $A_n$ of pairwise inequivalent elements has been listed. Then carry out the $n$-th stage of $\Gamma$ and enumerate a pair $(x,y)$. If $x\not\in F$ and there have been $l$ elements enumerated in $[x]_R$ by $\Gamma$, then we know that these $l$ elements form exactly $[x]_R$. In this situation list the smallest element of $[x]_R$, provided that it has not been listed before, and go on to the next stage. In the situation that the above conditions are not fulfilled, do nothing and go on to the next stage. From our assumption that $N_l=\infty$ it follows that this procedure would produce infinitely many elements, and from the construction they are pairwise inequivalent.
This finishes the proof for (1).

For (2), it suffices to show the conclusion for $k=2$. Toward a contradiction, assume there
is such a partial computable function $\rho$. Let $s$ be a computable function such 
that for any $e$ the ceer $R_{s(e)}$ is given by
$$ xR_{s(e)}y\Leftrightarrow x=y\lor (x=\rho(e,0)\!\downarrow\,\land 
y=\rho(e,1)\!\downarrow\,)\lor (x=\rho(e,1)\!\downarrow\,\land y=\rho(e,0)\!\downarrow\,),$$
for $x,y\in\NN$. 
It is obvious that $R_{s(e)}$ is 2-bounded for any $e$. Now by the fixed point theorem there is
some $e_0$ such that $W_{e_0}=W_{s(e_0)}$. For this $e_0$ we have that $R_{e_0}=R_{s(e_0)}$.
Since $R_{s(e)}$ is always 2-bounded, it follows that $R_{e_0}$ is 2-bounded and that
for all $x$, $\rho(e_0,x)\!\downarrow$. 
But then by the definition of $R_{s(e_0)}$ we have that $\rho(e_0,0)R_{s(e_0)}\rho(e_0,1)$,
which implies that $\rho(e_0,0)R_{e_0}\rho(e_0,1)$, which is a contradiction to our assumption that $\rho(e_0,\cdot)$ witnesses a reduction from $\omega$ to $R_{e_0}$.
\end{proof}

Fix $k>1$. Note that there is a computable enumeration of all $k$-bounded relations.
A canonical enumeration can be given as follows. For an index $e$, let $R_e$ be
the $e$-th ceer in the canonical enumeration of all ceers. We modify $R_e$ by
omitting any pair in its enumeration which would give rise to an equivalence
class of size $>k$. Denote this resulting $k$-bounded relation by $B^k_e$. Then
$\{ B^k_e\}_{e\in\omega}$ enumerates exactly all $k$-bounded relations.

For each $k>1$, define a ceer $B^k_{\infty}$ by
$$ \langle x,z\rangle B^k_{\infty}\langle y,z\rangle \Leftrightarrow
xB^k_zy, \mbox{ for $x,y,z\in \NN$.} $$
Then $B^k_{\infty}$ is $k$-bounded and is in fact universal for all $k$-bounded
relations, i.e., if $R$ is an arbitrary $k$-bounded relation, then $R\leq B^k_{\infty}$.

The following theorem establishes a hierarchy for bounded relations in analogy with
the hierarchy of finite dimensional relations.

\begin{theorem}[Bounded Hierarchy Theorem]
For any $k>0$, $B^k_\infty< B^{k+1}_\infty$.
\end{theorem}
\begin{proof}
It suffices to show that there is a $k+1$-bounded relation which is not essentially
$k$-bounded. For this, let $A$ be a non-computable c.e. set and define a ceer $R$ by
$$ \langle x,i\rangle R\langle x,j\rangle \Leftrightarrow i=j \lor (i,j\leq k \land x\in A),$$
for $x,i,j\in\NN$. Then $R$ is $k+1$-bounded. 

Now assume $R$ is essentially $k$-bounded, i.e., there are computable function $f$ and
$k$-bounded relation $S$ such that $R\leq S$ via $f$. Then we have
$$ x\in A\Leftrightarrow \exists i,j\leq k (i\neq j \land f(\langle x,i\rangle)=
f(\langle x,j\rangle)). $$
In fact, if $x\in A$, then $[\langle x,0\rangle]_R$ has $k+1$ elements; and since $f$ is
a reduction, it must map the $k+1$ elements into some single $S$-class, which only has
at most $k$-elements. On the other hand, if $x\not\in A$, then the elements $\langle x,i\rangle$
are pairwise $R$-inequivalent; again by that $f$ is a reduction, it follows that 
the statement on the right hand side cannot happen. Now the statement gives a computable
definition for $A$, a contradiction to our hypothesis that $A$ is non-computable.
\end{proof}

The idea of the above proof can be used to show that there are FC relations which are
not essentially bounded.

\begin{proposition}
There are FC relations which are not essentially bounded.
\end{proposition}
\begin{proof}
Again let $A$ be a non-computable c.e. set and define a ceer $R$ by a slight
modification to the previous definition, as follows. Let
$$ \langle x,i\rangle R\langle x,j\rangle \Leftrightarrow i=j \lor (i,j\leq x \land x\in A),$$
for $x,i,j\in\NN$. It is easy to see that $R$ is FC. 

Now assume $R$ is essentially $k$-bounded and let $f$ be a computable reduction function
witnessing this. Then by the same argument as before, we have that for $x>k$,
$$ x\in A\Leftrightarrow \exists i,j\leq x (i\neq j \land f(\langle x,i\rangle)=
f(\langle x,j\rangle)), $$
which gives again a computable definition for $A$, a contradiction.
\end{proof}

This proposition is also a corollary of Proposition 5.8 (1) and the fact that being
essentially bounded is $\Sigma_3$. But of course, the construction here makes
the proof more favorable than the abstract comparison of the complexities. Results
in {\S} 8 will show that if $A$ is hypersimple, then the FC relation $F_A$ defined
in Definition 5.4 is not essentially bounded. This provides another proof of
Proposition 6.4.

Similar
to Proposition 4.11, we can characterize the complexity of various index sets
related to boundedness.

\begin{proposition}
The following index sets are $\Sigma_3$-complete:
\begin{enumerate}
\item[(1)] $\{ e\,|\, R_e\mbox{ is essentially $n$-bounded}\}$, where $n>1$.
\item[(2)] $\{ e\,|\, R_e\mbox{ is essentially universal $n$-bounded}\}$, where $n>1$.
\item[(3)] $\{ e\,|\, R_e\mbox{ is essentially bounded}\}$.
\end{enumerate}
\end{proposition}

The proof is similar to that of Proposition 4.11. It is not clear whether $\{ e\,|\,
B^n_e\equiv B^n_\infty\}$ is $\Sigma_3$-complete for all $n>1$.

Next we further explore the notion of essential boundedness.

\begin{proposition}
Let $k\geq 3$. Let $R$ be a $k$-bounded relation such that every equivalence class
of $R$ is non-trivial. Then $R$ is essentially $\lfloor \frac{k}{2}\rfloor$-bounded.
In particular, if $k=3$, then $R$ is computable.
\end{proposition}
\begin{proof}
Let $\Gamma$ be an enumeration of the pairs of distinct elements in $R$. Without
loss of generality assume that at each stage $\Gamma$ enumerates one pair.
For each stage $s$, let $D_s$ be the set of all elements appeared
in some pair enumerated in $\Gamma$ before stage $s$ and let $E_s$ be
the equivalence relation on $D_s$ generated by $\Gamma$ before stage $s$. Then
we assume also that $E_{s+1}\neq E_s$ for every $s$. By our assumption that no
$R$-equivalence class is a singleton, $\bigcup_sD_s=\NN$.
We define a $\lfloor \frac{k}{2}\rfloor$-bounded relation $S$ and a reduction 
$f$ from $R$ to $S$ simultaneously by stages.

At stage $0$, let $S_0=\omega$ and $f_0$ be empty.
At stage $s+1$, suppose we have already defined $\lfloor \frac{k}{2}\rfloor$-bounded
$S_s$ on $D_s$ with $S_s\subseteq E_s$ and $f_s$ witnessing a reduction from $E_s$ to $S_s$. 
Let $(i,j)$ be the pair given by $\Gamma$ at stage $s$. There are four cases.
\begin{enumerate}
\item $D_{s+1}=D_s\cup \{ i,j\}$. In this case let $S_{s+1}$ be the trivial
extension of $S_s$ to the domain $D_{s+1}$ and let $f_{s+1}=f_s\cup \{ (i, \min(i,j)),
(j,\min(i,j))\}$.
\item $D_{s+1}=D_s\cup\{ i\}$. Let $C$ be the $E_s$-class of $j$. Define
$S_{s+1}$ be the trivial extension of $S_s$ to the domain $D_{s+1}$ and let 
$f_{s+1}=f_s\cup\{ (i, \min f_s(C))\}$.
\item $D_{s+1}=D_s\cup\{ j\}$. Similar to Case 2.
\item $D_{s+1}=D_s$. Let $C$ and $D$ be the $E_s$-classes of $i$ and $j$, respectively.
Then $C\cup D$ is an $E_{s+1}$-class. Define $S_{s+1}=S_s\cup \{ (\min f_s(C),\min f_s(D))\}$
and let $f_{s+1}=f_s$. 
\end{enumerate}
It is easy to see by induction that for each $S_s$-class $A$, there is some 
$E_s$-class $C$ such that $A\subseteq C$ and that the size of $A$ is less than
half of that of $C$. This guarantees that $S$ is $\lfloor \frac{k}{2}\rfloor$-bounded.
\end{proof}

A corollary to the proof is that all $3$-bounded relations are PC.

\begin{corollary}
All $3$-bounded relations are PC. In particular, $B^3_\infty\leq H$.
\end{corollary}
\begin{proof} Suppose $R$ is $3$-bounded. Use the same proof as above.
The relation $S$ constructed will be just $\omega$, since Case 4 never happens.
The only difference is that we do not know $\bigcup_sD_s=\NN$ now, corresponding
to the situation that $f$ may be only partial computable. But then $f$ witnesses
the definition of $R$ being PC. For the second assertion, just note that $H$ is
universal for PC relations.
\end{proof}

This corollary describes a part of the general relationship between bounded
relations and the halting equivalence and its jumps (see {\S} 8). In particular,
it is also true that $B^4_\infty\not\leq H$. We postpone the proof to {\S} 8, where
a much more general result is stated and proved.

To consider another generalization of the above results, we define spectra for ceers.
\begin{definition}
Let $R$ be ceer. Then the {\em spectrum} of $R$ is the set
$$ \spec(R)=\{ n\in\NN \,|\, n=|[x]_R| \mbox{ for some $x\in\NN$}\}.$$
\end{definition}

With this definition, Proposition 6.6 can be restated as follows:
If $R$ is $k$-bounded and $1\not\in\spec(R)$, then $R$ is essentially 
$\lfloor \frac{k}{2}\rfloor$-bounded.
The method above can be modified to show the following general results.

\begin{theorem}
Let $k>1$ and $R$ be a bounded relation.
\begin{enumerate}
\item[(1)] If $\spec(R)\subseteq [k,n]$ for some $n\geq k$, then
$R$ is essentially $\lfloor \frac{n}{k}\rfloor$-bounded.
\item[(2)] If $\spec(R)\subseteq \{ 1\}\cup [k,2k-1]$, then 
$R$ is PC.
\end{enumerate}
\end{theorem}
\begin{proof}
The proof of (1) is a modification of the
proof of Proposition 6.6, in which construction no action is taken until
an equivalence class of size $\geq k$ is generated.
The proof of (2) is the same as that of Corollary 6.7. 
\end{proof}

Nies \cite{Ni} considered 2-bounded relations and FC relations from a different
perspective. 

\begin{definition}
For equivalence relations $R_1$ and $R_2$ on $\NN$, let $R_1\lor R_2$ denote
the {\em join} of $R_1$ and $R_2$, which is the smallest equivalence relation
containing both $R_1$ and $R_2$.
\end{definition}

Nies proved a Join Theorem which states that
for any ceer $R$ there are 2-bounded relations $R_1$ and $R_2$ such that
$R=R_1\lor R_2$ and $R_1\cap R_2=\omega$.

\section{The saturation jump operator}

Throughout this section we will consider equivalence relations defined
on the set of all finite subsets of $\NN$, or in usual notation, on
$[\NN ]^{<\omega}$. There is a computable bijection between $[\NN ]^{<\omega}$ and $\NN$,
therefore we continue to regard the universe of our equivalence relations to be formally
$\NN$. However, it turns out to be less confusing to
work with the space $[\NN ]^{<\omega}$ for the equivalence relations
we are to define in this section.

\begin{definition}
Let $R$ be a ceer on $\NN$ and let $X\subseteq \NN$. The {\em
$R$-saturation} of $X$ is defined by
$$ [X]_R=\{ y\in\NN\,|\, \exists x\in X (xRy)\}. $$
\end{definition}

\begin{definition}
Let $R$ be a ceer on $\NN$. The {\em saturation jump} of $R$, 
denoted by $R^+$,
is an equivalence relation on $[\NN]^{<\omega}$ defined by
$$ xR^+y\Leftrightarrow  [x]_R=[y]_R. $$
For $n\in \NN$, the {\em $n$-th saturation jump} of $R$ is inductively defined by
$$\begin{array}{rcl}
R^{0+}&=& R \\
R^{1+}&=& R^+ \\
R^{(n+1)+}&=& (R^{n+})^+
\end{array}
$$
\end{definition}
Again, formally the universe of each $R^{n+}$ is $\NN$, but it is also natural to
indentify the elements as hereditarily finite sets of depth $n$. For example,
the universe of $R^{2+}$ could then be identified as the collection of finite sets
of finite sets.

The reader should be warned about the choice of the word ``jump" for the name of
the operator. Although it is obvious that $R\leq R^+$ always holds,
it is not true that for any ceer $R$, $R<R^+$. 
In fact, if $R$ is computable, then so is $R^+$. Below let us collect some basic facts
about the saturation jump operator.

\begin{proposition} Let $E$, $E_1$ and $E_2$ be arbitrary ceers and $n\in\NN$.
\begin{enumerate}
\item[(1)] If $E_1\leq E_2$, then $E^{n+}_1\leq E^{n+}_2$.
\item[(2)] If $E$ is FC, then so is $E^{n+}$.
\item[(3)] If $E$ is CC, then so is $E^{n+}$.
\item[(4)] If $E$ has the property that every pair of distinct $E$-classes are
computably separable, then so does $E^{n+}$.
\item[(5)] $E$ and $E^{n+}$ as sets have the same Turing degree.
\end{enumerate}
\end{proposition}
\begin{proof} Without loss of generality assume $n=1$.
(1)-(3) and (5) are obvious. For (4), let $\lnot xE^+y$ and assume without loss of generality
that $a\in x$ is such that $\forall b\in y$, $\lnot bEa$. For each $b\in y$, let 
$A_b$ be a computable set separating $[b]_E$ from $[a]_E$. Let $A=\bigcup_{b\in y} A_b$.
Then $A$ is computable. It follows that the set $\{ z\in[\NN]^{<\omega}\,|\, \exists u\in z
(u\not\in A)\}$ is computable and it separates $[x]_{E^+}$ from $[y]_{E^+}$.
\end{proof}

Next we consider essentially bounded relations and show that the saturation 
jump operator behaves properly on these relations.

\begin{theorem} If $E$ is essentially bounded and non-computable, then $E<E^+$.
\end{theorem}
\begin{proof}
Suppose $B$ is $k$-bounded and $f$ is computable witnessing that $E\leq B$.
Assume $g$ witnesses $E^+\leq E$. For notational simplicity we first consider $k=2$.
We derive a contradiction by showing that $E$ is computable. Given $x,y\in\NN$, note that 
$$\begin{array}{rcl}
xEy &\Leftrightarrow &
\mbox{ at least one pair of $\{ x\},\{ y\},\{ x,y\}$ are $E^+$-equivalent} \\
&\Leftrightarrow & \mbox{ all of $\{ x\},\{ y\},\{ x,y\}$ are $E^+$-equivalent} \\
&\Leftrightarrow & \mbox{ at least one pair of $g(\{ x\}),g(\{ y\}),g(\{ x,y\})$ are $E$-equivalent} \\
&\Leftrightarrow & \mbox{ all of $g(\{ x\}),g(\{ y\}),g(\{ x,y\})$ are $E$-equivalent} \\
&\Leftrightarrow & \mbox{ at least one pair of $f(g(\{ x\})), f(g(\{ y\})), f(g(\{ x,y\}))$
are $B$-equivalent} \\
&\Leftrightarrow & \mbox{ all of $f(g(\{ x\})), f(g(\{ y\})), f(g(\{ x,y\}))$
are $B$-equivalent} \\
&\Leftrightarrow &
\mbox{ at least one pair of $f(g(\{ x\})), f(g(\{ y\})), f(g(\{ x,y\}))$ are equal}
\end{array}
$$
where the last condition is computable. For general $k$, note that, by the same
argument as above,
$$\begin{array}{rcl}
xEy 
&\Leftrightarrow & \mbox{ at least one pair of 
$g(\{ x\}),g(\{ y\}),g(\{ x,y\})$ are $E$-equivalent} \\
&\Leftrightarrow & \mbox{ all of $g(\{ x\}),g(\{ y\}),g(\{ x,y\})$ are $E$-equivalent} \\
& & \mbox{ (letting $z_0=g(\{ x\}), z_1=g(\{ y\}), z_2=g(\{ x,y\})$} \\
&\Leftrightarrow & \mbox{ at least one pair of 
$g(\{ z_0\}), g(\{ z_1\}), g(\{ z_2\}), g(\{ z_0,z_1\}), g(\{ z_1,z_2\}),$} \\
& & \mbox{$g(\{ z_0,z_2\}), g(\{ z_0,z_1,z_2\})$ are $E$-equivalent} \\
&\Leftrightarrow & \mbox{ all of $g(\{ z_0\}), g(\{ z_1\}), g(\{ z_2\}), g(\{ z_0,z_1\}),
g(\{ z_1,z_2\}),$} \\
& & \mbox{$g(\{ z_0,z_2\}), g(\{ z_0,z_1,z_2\})$ are $E$-equivalent} \\
&\Leftrightarrow & \cdots\cdots
\end{array}
$$
When the iteration on the right hand side produces more than $k$ elements, say,
$w_0,\dots,w_m$, where $m\geq k-1$, we 
have that $xEy$ iff at least one pair of the elements $f(w_0),\dots,f(w_m)$ are equal.
This shows that $E$ is computable, a contradiction.
\end{proof}

Next we demonstrate that there are many ceers whose saturation jump is no more complicated.
For this we define a general operation $R^{\omega+}$ for ceer $R$ as follows.
First define by induction on $n$ a sequence of ceers $E_n$ by
$$\begin{array}{rcl}
\langle x,0\rangle E_0 \langle y,0\rangle & \Leftrightarrow & xRy \\
\langle x,i\rangle E_{n+1} \langle y,i\rangle & \Leftrightarrow & 
x=y \lor (i\leq n \land \langle x,i\rangle E_n\langle y,i\rangle) \lor
(i=n+1 \land xE_n^+y)
\end{array}
$$
Eventually define $R^{\omega+}$ by $$xR^{\omega+}y\Leftrightarrow \exists n (xE_ny).$$
It is then easy to see that $R^{n+}\leq R^{\omega+}$ for any $n\in\NN$ and that 
the properties of Proposition 7.3 remain true for $n=\omega$.

\begin{proposition}
For any ceer $R$, $(R^{\omega+})^+\equiv R^{\omega+}$. In particular, there are FC relations
$R$ such that $R^+\equiv R$.
\end{proposition}
\begin{proof}
Let $E=R^{\omega+}$. For each $x$ let $k_x$ be the largest $i$ such that 
$\langle y,i\rangle \in x$ for some $y$. Define a reduction function $f$ of $E^+$ to $E$ by
$f(x)=\langle x,k_x+1\rangle$. This $f$ works, since if $uE^+v$ then $k_u=k_v$ and hence 
$f(u)E_{k_u+1}f(v)$. On the other hand if $\lnot (uE^+v)$ then either $k_u\neq k_v$ or
else $\lnot(uE_{k_u}v)$, in either case $\lnot (f(u)Ef(v))$. 
\end{proof}

The following fact shows that the saturation jump operation is not one-one
on bi-reducibility degrees of ceers.
\begin{proposition}
There are ceers $R_1$ and $R_2$ such that $R_1<R_2$ but $R_1^+\equiv R_2^+$.
\end{proposition}
\begin{proof}
Let $R_1$ and $R_2$ be respectively defined by
$$ \langle x,0\rangle R_1 \langle x,1\rangle \Leftrightarrow x\in K $$
and 
$$ \langle x,0\rangle R_2 \langle x,1\rangle R_2 \langle x,2\rangle \Leftrightarrow x\in K.$$
Then $R_1\leq R_2$, $R_1$ is $2$-bounded and $R_2$ is not essentially $2$-bounded by the
proof of Theorem 6.3, hence
$R_2\not\leq R_1$.

To see that $R_1^+\equiv R_2^+$, it is enough to demonstrate a reduction from $R_2^+$
to $R_1^+$. For this let $g_1$ and $g_2$ be two computable functions with the property
that $g_1(x)\neq g_2(y)$, $\forall x,y$ and $x\in K\Leftrightarrow g_1(x)\in K
\Leftrightarrow g_2(x)\in K$. Such functions are easy to construct.
Now define a computable function $f$ such that
$$\begin{array}{rcl}
f(\langle x,0\rangle)&=&\{ \langle g_1(x),0\rangle, \langle g_2(x),0\rangle\}, \\
f(\langle x,1\rangle)&=&\{ \langle g_1(x),0\rangle, \langle g_2(x),1\rangle\}, \mbox{ and} \\
f(\langle x,2\rangle)&=&\{ \langle g_1(x),1\rangle, \langle g_2(x),1\rangle\}.
\end{array}
$$
It is easy to see that $f$ induces a required reduction.
\end{proof}

Finally some curious examples of the effect of the saturation jump operator.

\begin{theorem}
Let $H$ be the halting equivalence. Then $H<H^+$.
\end{theorem}
\begin{proof}
We first show that if $H^+\leq H$, then $H^+$ is PC. Suppose a computable function
$f$ witnesses $H^+\leq H$ and a partial computable function $\varphi$ witnesses
$H$ is PC. We claim that $\varphi\circ f$ witnesses that $H^+$ is PC. For this
it is enough to check that for any $x,y\in[\NN]^{<\omega}$ with $\varphi(f(x))\!\uparrow\,$
and $f(x)=f(y)$, we have $x=y$. By the definition of $H$, if $\varphi(f(x))\!\uparrow\,$,
then $[f(x)]_H$ is a singleton. Since $f$ is a reduction, it follows that $[x]_{H^+}$
is computable. Therefore $x$, as a finite subset of $\NN$, consists only of elements
$a$ with $\varphi_a(a)\!\uparrow\,$. It follows then $[x]_{H^+}$ is a singleton, thus
$f(x)=f(y)$ implies that $x=y$.

Next we show that $H^+$ is not PC. Assume a partial computable function $\psi$ witnesses 
that $H^+$ is PC. Then if $x\in [\NN ]^{<\omega}$ has more than one element in its $H^+$-class,
we have $\psi(x)\!\downarrow\,$. Now pick an arbitrary element $a\in K_0$, i.e., with 
$\varphi_a(a)\!\downarrow\, =0$. Let $b=\psi(\{ a\})$. Consider the function $f$ defined
by $f(i)=\psi(\{ a, i\})$, for $i\in\NN$. Then $f$ is a total computable function.
Moreover, $i\in K_0$ iff $f(i)=b$. Since $K_0$ is not computable, this is a contradiction.
\end{proof}

\begin{theorem}
For any $n\in\NN$, $H^{n+}$ is not universal, i.e., $H^{n+}<R_{\infty}$.
\end{theorem}
\begin{proof}
Let $R$ be a precomplete ceer. Then every pair of distinct $R$-equivalence classes 
are computably inseparable (see, e.g. \cite{BS}). 
We show that $R\not\leq H^{n+}$ for any $n\in \NN$. For $n=0$ this
was done in Proposition 3.9. We assume $n\geq 1$. We use the union operation in
the sense of set theory, i.e., for a set $x$, let
$$ \bigcup x =\{ z\,|\, \exists y\in x (z\in y)\}. $$
For $k\geq 1$, define inductively $\bigcup^{k+1}x=\bigcup (\bigcup^kx)$.
Now fix $n\geq 1$ and assume that $f$ is a computable function witnessing $R\leq H^{n+}$.
For each $x\in\NN$, let $k(x)=\bigcup^{n-1}f(x)\cap\overline{K}$, then each $k(x)$ is a
finite subset of $\overline{K}$. We claim that $k$ is a constant function. Since
otherwise, there are $x\neq y$ such that $k(x)\neq k(y)$. Without loss of generality
assume $a\in k(x)\setminus k(y)$. Then $f(x)$ and $f(y)$ are not equivalent under $H^{n+}$,
therefore $\lnot xRy$. It follows that $[f(x)]_{H^{n+}}$ and $[f(y)]_{H^{n+}}$ are
computably inseparable. But this is not so, since the computable set 
$\{ z\,|\, a\in\bigcup^{n-1}z\}$ separates them, a contradiction. 

Now suppose the function $k$ takes the constant value $\{ c_1,\dots,c_m\}$.
Then for each $x\in\NN$, $\bigcup^{n-1}f(x)\setminus\{ c_1,\dots,c_m\}$ consists
only of elements of $K$. It follows that there is a computable function $g$ witnessing
that $R\leq \omega$, a contradiction.
\end{proof}

A similar argument, with a bit more notation, shows the same for $n=\omega$, 
but we omit it now since we will demonstrate a stronger result later.

\begin{theorem}
Let $E$ be a finite dimensional ceer. Then for any $n\in\NN$, $H\not\leq E^{n+}$.
\end{theorem}
\begin{proof} This is a similar argument as in the preceding proof.
Without loss of generality assume $E=R_{K_0,\dots,K_{m-1}}$ is universal $m$-dimensional
and let $R=E^{n+}$. Assume toward a contradiction that $H\leq R$ via $f$. For each $x\in\NN$,
let $A_x=\bigcup^{n-1}f(x)\cap (K_0\cup\dots\cup K_{m-1})$ and 
$B_x=\bigcup^{n-1}f(x)\setminus A_x$.
We claim that $B_x$ is a constant finite set for all $x\in K$. This is because, if $B_x\neq B_y$
for some $x\neq y\in K$, then $[f(x)]_R$ and $[f(y)]_R$ can be computably separated, hence
$[x]_H$ and $[y]_H$ can be computably separated, which is false. Now it follows that there is
an infinite set $S\subseteq K$ such that for any $x,y\in S$, $x\neq y$, $[A_x]_E\neq [A_y]_E$.
But this is impossible since there are at most $2^m$ many distinct $[A]_E$ with $A\subseteq
K_0\cup\dots\cup K_{m-1}$.
\end{proof}

\begin{theorem} Let $A$ be hypersimple. For any $n\in\NN$, $\omega\not\leq F_A^{n+}$.
\end{theorem}
\begin{proof}
We have seen in Proposition 5.6 that the conclusion is true for $n=0$.
Suppose $n>0$ and assume $\omega\leq F_A^{n+}$ via $f$. We derive a contradiction
by defining a reduction $g$ from $\omega$ to $F_A$. Let $g(0)=\sup\bigcup^{n-1}f(0)$.
In general we will define each $g(m)$ as $\sup\bigcup^{n-1}\bigcup_{i\leq k}f(i)$ for some $k$.
Suppose $g(m)$ is defined.

Let $s(z,n)$ be defined by the recursion $s(z,0)=z$ and $s(z,n+1)=2^{s(z,n)}$. Then
$s(z,n)$ dominates the number of $F_A^{n+}$-equivalence classes in which
there is an element $x$ so that $\sup\bigcup^{n-1}x\leq z$. 

Now let 
$g(m+1)=\sup\bigcup^{n-1}\bigcup_{i\leq s(g(m),n)+1}f(i)$. By the property of $s$
mentioned above and the assumption on $f$, there is some element $x$ with $g(m)<x\leq g(m+1)$
such that $\lnot (g(m)F_Ax)$. This implies that $\lnot (g(m+1)F_Ag(m))$.
\end{proof}

\section{The halting jump operator}

In this section we introduce another jump operator which reminds us much of
the halting jump in computability theory.

\begin{definition}
Let $E$ be an arbitrary ceer. The {\em halting jump} of $E$, denoted by $E'$,
is the ceer defined by
$$ xE'y\Leftrightarrow x=y\lor \varphi_x(x)\!\downarrow\, E\varphi_y(y)\!\downarrow\, $$
for $x,y\in\NN$. For any $n\in\NN$, the {\em $n$-th halting jump} of $E$, denoted
by $E^{(n)}$, is inductively defined by
$$ \begin{array}{rcl}
E^{(0)}&=& E \\
E^{(1)}&=& E' \\
E^{(n+1)}&=& (E^{(n)})'
\end{array}
$$
\end{definition}

\begin{definition}
Let $E$ be an arbitrary ceer. A ceer $R$ is {\em partially classifiable by $E$},
denoted by {\em $\mbox{PC}^E$}, if there is a partial computable function $\psi$ such
that, for any $x,y\in\NN$,
$$ xRy\Leftrightarrow x=y \lor \psi(x)\!\downarrow\, E\psi(y)\!\downarrow\, .$$
\end{definition}

Then we have the following basic properties.

\begin{proposition} Let $E$ be an arbitrary ceer.
\begin{enumerate}
\item[(1)] $E\leq E'$.
\item[(2)] $E'$ is universal for $\mbox{PC}^E$ relations.
\item[(3)] For $n\geq 1$, $\mbox{id}(n)'\equiv R_{K_0,\dots,K_{n-1}}$.
\item[(4)] $\omega '=H$.
\item[(5)] Let $R$ be a ceer without computable classes. If $R\leq E'$ then $R\leq E$
\end{enumerate}
\end{proposition}
\begin{proof}
(1), (3) and (4) are obvious. (2) uses the same proof as that of Proposition 3.8.
(5) is proved in the same way as Proposition 3.9.
\end{proof} 

Proposition 8.3 (5) implies that the halting jump of a non-universal ceer is non-universal.
This feature of the halting jump operator is also a corollary of the following
interesting result.

\begin{theorem}[The Halting Jump Theorem]
For ceers $E_1$ and $E_2$, $E_1\leq E_2$ iff $E'_1\leq E'_2$.
\end{theorem}
\begin{proof}
First suppose $E_1\leq E_2$. Let $f$ be a computable function such that
$xE_1y\Leftrightarrow f(x)E_2f(y)$, for any $x,y\in\NN$. First
define a one-one computable function $g$ so that for any $x,z\in\NN$,
$\varphi_{g(x)}(z)=f(\varphi_x(x))$. Then $\varphi_{g(x)}(g(x))\!\downarrow\,$ iff
$f(\varphi_x(x))\!\downarrow\,$, and when this happens, they take the same value.
Therefore, since $g$ is one-one, we have
$$\begin{array}{rcl}
xE'_1y &\Leftrightarrow & x=y \lor \varphi_x(x)\!\downarrow\, E_1 \varphi_y(y)\!\downarrow\, \\
&\Leftrightarrow & x=y \lor f(\varphi_x(x))\!\downarrow\, E_2 f(\varphi_y(y))\!\downarrow\, \\
&\Leftrightarrow & g(x)=g(y) \lor \varphi_{g(x)}(g(x))\!\downarrow\, 
E_2 \varphi_{g(y)}(g(y))\!\downarrow\, \\
&\Leftrightarrow & g(x)E'_2 g(y)
\end{array}
$$
This shows that $E'_1\leq E'_2$ via $g$.

Next suppose $E'_1\leq E'_2$ via $f$. 
Then we have that for any $x\in K$, $f(x)\in K$. This is because, if $f(x)\not\in K$,
then $[f(x)]_{E'_2}$ is a singleton. Since $f$ is a reduction, it follows that
$[x]_{E'_1}$ is computable; but this can only happen when $x\not\in K$ (if $x\in K$,
then it is easy to see that $[x]_{E'_1}$ would be $m$-complete). 

Now let us fix a one-one computable function $s$ such that $\varphi_{s(x)}(s(x))=x$.
Then define a computable function $g$ by $g(x)=\varphi_{f(s(x))}(f(s(x)))$.
$g$ is total since for any $x$, $s(x)\in K$ by definition, and therefore $f(s(x))\in K$
by the claim in the preceding paragraph. Finally,
$$ \begin{array}{rcl}
xE_1y & \Leftrightarrow & s(x)E'_1s(y) \mbox{ (by definition of $s$)} \\
& \Leftrightarrow & f(s(x))E'_2f(s(y)) \\
& \Leftrightarrow & g(x)E_2g(y) \mbox{ (by definition of $E'_2$ and $g$)}
\end{array}
$$
Therefore $E_1\leq E_2$.
\end{proof}

\begin{corollary}$~$\
\begin{enumerate}
\item[(1)] For ceers $E_1$ and $E_2$, $E_1\equiv E_2$ iff $E'_1\equiv E'_2$.
\item[(2)] If $E$ is non-universal, then $E'$ is non-universal.
\item[(3)] For any ceer $E$, exactly one of the following holds:
\begin{enumerate}
\item[i)] for any $n\in\NN$, $E^{(n)}\equiv E$;
\item[ii)] for any $n\in\NN$, $E^{(n)}<E^{(n+1)}$.
\end{enumerate}
\item[(4)] For any $n\in \NN$, $H^{(n)}<H^{(n+1)}$.
\item[(5)] For any $n\in\NN$, the sets $\{ e\,|\, R_e\leq H^{(n)}\}$,
$\{ e\,|\, R_e\geq H^{(n)}\}$ and $\{ e\,|\, R_e\equiv H^{(n)}\}$ are all
$\Sigma_3$-complete.
\item[(6)] The set $\{ e\,|\, R_e\mbox{ is PC}\}$ is $\Sigma_3$-complete.
\end{enumerate}
\end{corollary}
\begin{proof}
(1) and (2) are immediate. 
For (3), note that $E\leq E'\leq\dots\leq E^{(n)}\leq E^{(n+1)}\leq\dots$.
If ii) fails, then for some $n\in\NN$, $E^{(n+1)}\leq E^{(n)}$; by applying Theorem
8.4 $n$ times, we have $E'\leq E$, therefore $E\equiv E'$, which implies that
$E^{(n)}\equiv E$ for any $n\in\NN$. (4) follows immediately from (3) since $\omega<H=\omega'$.
(5) follows immediately from Propositions 4.7-4.10. For (6) a straightforward computation
shows that the set is $\Sigma_3$. It is then easy to check that if $W_e$ is computable,
then $R_{W_e}$ is PC. The Halting jump theorem implies that if $W_e$ is not computable,
then $R_{W_e}\not\leq\omega'=H$, hence in particular it is not PC. Thus we have a reduction
of Comp to $\{ e\,|\, R_e\mbox{ is PC}\}$.
\end{proof}

The next few theorems deal with bounded relations and their halting jumps.
\begin{theorem}
Let $n\in\NN$ and $R$ be an $n$-bounded relation. Then there is an 
$\lfloor \frac{n}{2}\rfloor$-bounded relation $S$ such that $R\leq S'$. 
Equivalently, $B^n_\infty\leq (B^{\lfloor \frac{n}{2}\rfloor}_\infty)'$.
\end{theorem}
\begin{proof}
It is enough to define an $\lfloor \frac{n}{2}\rfloor$-bounded relation $S$
so that $R\in \mbox{PC}^S$, i.e., for some partial computable function $\psi$,
$$ xRy\Leftrightarrow x=y \lor \psi(x)\!\downarrow\, S\psi(y)\!\downarrow\,.$$
We fix an enumeration of the pairs of distinct elements in $R$ and 
define $S$ and $\psi$ simultaneously.
At the stage $s$ of the enumeration let the pair $(i,j)$ (say $i<j$) be enumerated 
for $R$ and suppose $S_s$ and $\psi_s$ are defined. Consider four cases:
\begin{enumerate}
\item $i,j\not\in\dom\psi_s$. In this case extend
$\psi_s$ by letting $\psi_{s+1}(i)=\psi_{s+1}(j)=i$ and let $S_{s+1}=S_s$.
\item $i\in\dom\psi_s$ but $j\not\in\dom\psi_s$. In this case extend
$\psi_s$ by letting $\psi_{s+1}(j)=\psi_s(i)$ and let $S_{s+1}=S_s$.
\item $i\not\in\dom\psi_s$ but $j\in\dom\psi_s$. In this case similarly
extend $\psi_s$ by letting $\psi_{s+1}(i)=\psi_s(j)$ and let $S_{s+1}=S_s$.
\item $i,j\in\dom\psi_s$. If in addition $\psi_s(i)=\psi_s(j)$ then do
nothing; just let $\psi_{s+1}=\psi_s$ and $S_{s+1}=S_s$. Otherwise, $\psi_s(i)\neq
\psi_s(j)$. Then let $\psi_{s+1}=\psi_s$ and $S_{s+1}=S_s\cup\{ (\psi_s(i),\psi_s(j))\}$.
\end{enumerate}
This finishes our construction of $S$ and $\psi$. It follows easily from the
construction that $R$ is $\mbox{PC}^S$. Note that for any $x\in\dom\psi$, $\psi(x)Rx$;
in addition, there is $y\neq x$ with $\psi(y)=\psi(x)$. Also note that $S$ admits new 
pairs only in Case 4 above, which implies that $S_s\subseteq R_s$ and that
every $S_s$-class contains at most half many elements as the corresponding $R_s$-class. 
This implies that $S$ is $\lfloor \frac{n}{2}\rfloor$-bounded.
\end{proof}

\begin{corollary}
For any $n\in\NN$ and $(2^{n+1}-1)$-bounded relation $R$, $R\leq \omega^{(n)}$. Equivalently, $B^{2^{n+1}-1}_\infty\leq\omega^{(n)}$.
\end{corollary}
\begin{proof}
By induction on $n$. The case $n=0$ is obvious. Suppose $R$ is $(2^{n+2}-1)$-bounded.
Then by the previous theorem there is $(2^{n+1}-1)$-bounded relation $S$ such that $R\leq S'$.
By the inductive hypothesis $S\leq \omega^{(n)}$. Therefore $R\leq S'\leq \omega^{(n+1)}$.
\end{proof}

It might be worth noticing that Proposition 6.6 is a corollary of Theorem 8.6, and that Corollary 6.7 is a special case of Corollary 8.7. In fact, if $R$ is $n$-bounded and
every equivalence class of $R$ is non-trivial, then the proof of Theorem 8.6 produces
a {\em total} function $\psi$ witnessing the reduction from $R$ to some $\lfloor \frac{n}{2}\rfloor$-bounded ceer.

Next we show that the above estimates are in some sense optimal. 
In dealing with iterated halting jumps 
it is helpful to introduce the following 
notation. Let $\kappa$ denote the partial computable function
$\kappa(x)=\varphi_x(x)$, for any $x\in\NN$, and for any $n\in\NN$ let $\kappa^n$
be the $n$-th iterate of $\kappa$.
Then it is easy to see that, for any $x,y\in\NN$,
$$ x\omega^{(n)}y\Leftrightarrow \exists i\leq n (\kappa^i(x)\!\downarrow\, =\kappa^i(y)
\!\downarrow\,).$$

We will use the following lemma.
\begin{lemma} Let $C$ be a computable subset of $K$. Then there is a total computable
function $s(x)$ such that $\kappa(s(x))\!\downarrow\,=x$ and $s(x)\not\in C$, for any $x\in\NN$.
\end{lemma}
\begin{proof}
For $i\in\NN$ let $K_i=\{ x\,|\,\kappa(x)\!\downarrow\,=i\}$. Then each $K_i$ is non-computable.
Note that for any $i\in\NN$, $K_i\not\subseteq C$, since otherwise $K_i$ would have a
computable definition from the computability of $C$. This shows that for any $x\in\NN$,
the set $\{ y\,|\, \kappa(y)\!\downarrow\,=x \land y\not\in C\}$ is nonempty. Thus we may
construct a total computable $s$ so that for each $x$ it returns an element from this set.
\end{proof}

\begin{theorem} For any $n\in\NN$, $B^{2^{n+1}}_\infty\not\leq \omega^{(n)}$.
\end{theorem}
\begin{proof}
We define a sequence of bounded relations $E_n$ by induction on $n$. Let $E_0$ be
defined by
$$ \langle x,0\rangle E_0\langle x,1\rangle\Leftrightarrow \kappa(x)\!\downarrow\,. $$
Given $E_n$, define $E_{n+1}$ by
$$ \begin{array}{rcl}
\langle x,i\rangle E_{n+1}\langle x,j\rangle & \Leftrightarrow &
\langle x,i\rangle E_n\langle x,j\rangle \lor 
\langle x,i-2^{n+1}\rangle E_n\langle x,j-2^{n+1}\rangle \lor \\
& & (i,j<2^{n+2} \land \kappa^{n+2}(x)\!\downarrow\,).
\end{array}
$$
It is easy to see that each $E_n$ is $2^{n+1}$-bounded. We verify below that
$E_n\not\leq\omega^{(n)}$ for any $n\in\NN$, therefore establishing the theorem.

For $n=0$ the conclusion is obvious. Now fix an arbitrary $n>0$ and assume toward
a contradiction that $E_n\leq \omega^{(n)}$ via some function $f$. We define by induction
on $k\leq n$ computable sets $C_k$ and computable functions $s_k$ as follows.
Let $C_0=\{ x\,|\, \exists l<2^n (f(\langle x,2l\rangle)=f(\langle x,2l+1\rangle))\}$.
Then $C_0$ is computable; also $C_0\subseteq K$ since for any $x\in C_0$ there is some
$l<2^n$ such that $\langle x,2l\rangle E_n\langle x,2l+1\rangle$, whereas if $x\not\in K$
then $\lnot(\langle x,i\rangle E_n\langle x,j\rangle)$ whenever $i\neq j$.
Let $s_0$ be obtained from the preceding lemma for $C_0$. Then for any $l<2^n$ and any
$x$, we have that $f(\langle s_0(x),2l\rangle)\neq f(\langle s_0(x),2l+1\rangle)$;
but since $s_0(x)\in K$ and $\langle s_0(x),2l\rangle E_n \langle s_0(x),2l+1\rangle$
for any $l<2^n$,
it follows that for any $i<2^{n+1}$, $\kappa(f(\langle s_0(x),i\rangle))\!\downarrow\,$.
This finishes the definition for the base step. 

For the inductive step suppose we have defined computable set 
$C_k\subseteq K$ and computable function $s_k$ for $0\leq k< n$ with the following
properties for any $x\in\NN$:
\begin{enumerate}
\item[(a)] $\kappa^{k+1}(s_k(x))\!\downarrow\,=x$;
\item[(b)] For any $i<2^{n+1}$, $\kappa^{k+1}(f(\langle s_k(x),i\rangle))\!\downarrow\,$.
\end{enumerate}

Define $C_{k+1}$ to be the set
$$\{ x\,|\,\exists l<2^{n-k-1} (f(\langle s_k(x),l2^{k+2}\rangle)\,\omega^{(k+1)}\,
f(\langle s_k(x),l2^{k+2}+2^{k+1}\rangle))\}.$$ 
It follows from property (b) that $C_{k+1}$ is computable. To see that
$C_{k+1}\subseteq K$, it is enough to show that for $x\in C_{k+1}$,
$\kappa^{k+2}(s_k(x))\!\downarrow\,$. Fix $x\in C_{k+1}$ and a witness $l<2^{n-k-1}$.
Note that for $i=l2^{k+2}$ and $j=l2^{k+2}+2^{k+1}$, we have 
$f(\langle s_k(x),i\rangle)\omega^{(k+1)}f(\langle s_k(x),j\rangle)$ 
by the definition of $C_{k+1}$. It follows that 
$f(\langle s_k(x),i\rangle)\omega^{(n)}f(\langle s_k(x),j\rangle)$
and therefore $\langle s_k(x),i\rangle E_n\langle s_k(x),j\rangle$ since $f$ witnesses
$E_n\leq\omega^{(n)}$. But by the definition of $E_n$, it is easy to see that
if $\langle y,i\rangle E_n\langle y,j\rangle$ and $j-i=2^{k+1}$, then in fact
$\langle y,i\rangle E_{k+1}\langle y,j\rangle$. Thus we actually have that
$\langle s_k(x),i\rangle E_{k+1}\langle s_k(x), j\rangle$.
On the other hand, by the definition of $E_k$ and the fact that $E_k$ is
$2^{k+1}$-bounded, $\lnot (\langle s_k(x),i\rangle E_k \langle s_k(x),j\rangle)$
and $\lnot (\langle s_k(x),i-2^{k+1}\rangle E_k \langle s_k(x),j-2^{k+1}\rangle)$. 
It follows from the definition of $E_{k+1}$ that we must have $\kappa^{k+2}(s_k(x))\!\downarrow\,$.

Let $s'$ be obtained from the preceding lemma for $C_{k+1}$ and let 
$s_{k+1}(x)=s_k(s'(x))$. We verify that properties (a) and (b) are preserved by this
definition. (a) is immediate, since $$\kappa^{k+2}(s_{k+1}(x))\!\downarrow\,=\kappa(\kappa^{k+1}(s_k(s'(x))))\!\downarrow\,
=\kappa(s'(x))\!\downarrow\,=x. 
$$
For (b), first note that for any $i<2^{n+1}$, 
$\kappa^{k+1}(f(\langle s_{k+1}(x),i\rangle))\!\downarrow\,$ by the definition of $s_{k+1}$.
Thus if $k+1=n$ there is nothing to check. Suppose $k+1<n$. For any $x$, $s'(x)\not\in C_{k+1}$,
therefore, for any $l<2^{n-k-1}$, $p=l2^{k+2}$ and $q=p+2^{k+1}$, we have that  
$\lnot (f(\langle s_{k+1}(x),p\rangle)\,\omega^{(k+1)}\,f(\langle s_{k+1}(x),q\rangle))$.
By the inductive definition of $E_{k+1}$ this implies that for any $i<2^{n+1}$
there is $j<2^{n+1}$, $j\neq i$, such that 
$\lnot (f(\langle s_{k+1}(x),i\rangle)\,\omega^{(k+1)}\,f(\langle s_{k+1}(x),j\rangle))$.
By property (a) for $s_{k+1}$, we know that $\langle s_{k+1}(x),i\rangle E_n
\langle s_{k+1}(x),j\rangle$, for all $i,j<2^{n+1}$. Therefore
$f(\langle s_{k+1}(x),i\rangle)\,\omega^{(n)}\,f(\langle s_{k+1}(x),j\rangle)$,
for all $i,j<2^{n+1}$. Thus $\kappa^{k+2}(f(\langle s_{k+1}(x),i\rangle)\!\downarrow\,$
for any $i<2^{n+1}$.

In effect we have defined $C_k$, $s_k$ which satisfy properties (a) and (b) for $k<n$,
and $C_n$, $s_n$ with property (a) only. We have
that $\langle s_n(x),0\rangle E_n \langle s_n(x),2^n\rangle$ by the property (a) for $s_n$. 
However, for any $x$, we also have $\lnot(f(\langle s_n(x),0\rangle)\,\omega^{(n)}\,$
$f(\langle s_n(x),2^n\rangle))$. This contradicts the assumption that $f$ is a reduction
from $E_n$ to $\omega^{(n)}$.
\end{proof}

We have established in Corollary 8.7 that all bounded relations are reducible to some $\omega^{(n)}$.
Next we show that FC relations are not necessarily so. 
\begin{theorem}
Let $A$ be hypersimple. Then for any $n\in\NN$, $F_A\not\leq \omega^{(n)}$.
\end{theorem}
\begin{proof}
By induction on $n$. For $n=0$ the conclusion is trivial. Suppose now $F_A\not\leq\omega^{(n)}$
but $F_A\leq\omega^{(n+1)}$ via some computable function $f$. Therefore for any $x,y\in\NN$,
$$ xF_Ay\Leftrightarrow \exists i\leq n+1 (\kappa^i(f(x))\!\downarrow\, =\kappa^i(f(y))
\!\downarrow\,).$$
Note that there are infinitely many pairs $(x,y)$ such that $xF_Ay$ but
$\kappa^n(f(x))\!\downarrow\,\neq\kappa^n(f(y))\!\downarrow\,$. This is because, if it were
not the case, then one could modify $f$ by changing its values at finitely many points
so as to get that $F_A\leq\omega^{(n)}$, contradicting our hypothesis. Also note that
the (infinite) set of all pairs with the above property is c.e. Let us fix a non-repeating
computable enumeration of this set and denote its $k$-th element by $(x_k,y_k)$.

Now we define a computable function $g$ as follows. Let $g(0)=x_0$. For $i>0$ let
$g(i)$ be the first element of the sequence $x_0,y_0,x_1,y_1,\dots$ such that
for any $j<i$, $\kappa^{n+1}(f(g(j))\neq\kappa^{n+1}(f(g(i))$. Note that
the definition makes sense because for each element $z$ in the above sequence
$\kappa^{n+1}(f(z))\!\downarrow\,$, and because $F_A$ is FC. 

Thus $g$ is a total computable function with the property that if $i\neq j$,
then $\kappa^{n+1}(f(g(i))\!\downarrow\,\neq\kappa^{n+1}(f(g(j))\!\downarrow\,$,
hence $\lnot(g(i)F_Ag(j))$. This shows that $\omega\leq F_A$ via $g$, a contradiction
to Proposition 5.6.
\end{proof}

\begin{corollary} If $A$ is hypersimple, then $F_A$ is not essentially bounded. 
\end{corollary}
\begin{proof} This follows from the preceding theorem and Corollary 8.7.
\end{proof}

Corollary 8.5 (6) can also be seen as an immediate corollary of Theorem 8.10 by
using the theorem of Dekker we quoted earlier in the proof of Lemma 4.9.

From the alternative definition of $\omega^{(n)}$ using partial computable functions
$\kappa$ and its iterates we can define a natural limit of the 
$\omega^{(n)}$'s when $n$ approaches
infinity. Let $\omega^{(\omega)}$ be the ceer defined by
$$ x\omega^{(\omega)}y\Leftrightarrow \exists i 
(\kappa^i(x)\!\downarrow\,=\kappa^i(y)\!\downarrow\,)$$
for any $x,y\in\NN$. One may wonder if $F_A\leq \omega^{(\omega)}$ for hypersimple
$A$. We show next that this is indeed the case. In fact, $\omega^{(\omega)}$ is
universal for all ceers.

We need the following lemma.

\begin{lemma}
For any partial computable function $\psi$ there is a one-one total computable function
$v$ such that
$$ \kappa(v(x))=v(\psi(x)), \mbox{ for any $x\in\NN$.}$$
Moreover, an index for $v$ can be obtained effectively from an index of $\psi$.
\end{lemma}
\begin{proof}
Let $s(e,i)$ be a one-one total computable function such that
$$ \varphi_{s(e,\langle y,x\rangle)}(u)=\varphi_e(\varphi_y(x)). $$
Let $t(e)$ be a computable function such that $\varphi_{t(e)}(i)=s(e,i)$.
By the fixed point theorem there is $e_0$ such that $\varphi_{t(e_0)}=\varphi_{e_0}$.
Then
$$\varphi_{s(e_0,\langle y,x\rangle)}(u)=\varphi_{e_0}(\varphi_y(x))
=\varphi_{t(e_0)}(\varphi_y(x))=s(e_0,\varphi_y(x)).$$
Let $r(y)$ be a computable function such that $\varphi_{r(y)}(x)=\langle y,\psi(x)\rangle$.
By the fixed point theorem again, there is $y_0$ such that $\varphi_{r(y_0)}=\varphi_{y_0}$.
Then $\varphi_{y_0}(x)=\varphi_{r(y_0)}(x)=\langle y_0,\psi(x)\rangle.$
Therefore
$$ \varphi_{s(e_0,\langle y_0,x\rangle)}(u)=s(e_0,\varphi_{y_0}(x))=s(e_0,\langle y_0,\psi(x)
\rangle.$$
Letting $v(x)=s(e_0,\langle y_0,x\rangle)$, it is easy to check that $v$ is as required.
\end{proof}
\begin{theorem}
$\omega^{(\omega)}$ is universal for all ceers.
\end{theorem}
\begin{proof}
Let $R$ be an arbitrary ceer. Fixing a computable enumeration of $R$, let $R_s$ be
the equivalence relation on $\NN$ established by stage $s$ of the enumeration, i.e.,
$R_s$ is generated by the finitely many pairs enumerated up to stage $s$. Let
$p_i$ be the $i+1$-th prime number. Define a total computable function $\psi$ by
$\psi(p_i^s)=p_j^{s+1}$ if $j\leq i$ is the least element such that $jR_si$, and
$\psi(n)=0$ for any $n$ not a power of a prime. Let $v$ be the one-one total computable function given by the
preceeding lemma. Let $f(i)=v(p_i)$. We claim that $f$ is a reduction function from
$R$ to $\omega^{(\omega)}$.

It is easy to see that, for any $i,j\in\NN$, 
$$iRj\Leftrightarrow \exists s (\psi^s(p_i)\!\downarrow\, =\psi^s(p_j)\!\downarrow\,).$$
By induction on $s$ one can also see that 
$\kappa^s(f(i))=\kappa^s(v(p_i))=v(\psi^s(p_i))$.
Since $v$ is one-one, we have that $\psi^s(p_i)\!\downarrow\,=\psi^s(p_j)\!\downarrow\,$ iff
$\kappa^s(f(i))\!\downarrow\,=\kappa^s(f(j))\!\downarrow\,$. 
\end{proof}

We do not know if there are non-universal ceers $R$ with $R'\equiv R$. But the following
proposition suggests that it is probably unlikely.
\begin{proposition}
If $R'\equiv R$ then $\omega^{(n)}<R$ for any $n\in\NN$.
\end{proposition}
\begin{proof}
Since $1\leq R$ and $\omega<1'=R_K$, we have that 
$\omega^{(n)}<1^{(n+1)}\leq R^{(n+1)}$. But if $R'\equiv R$ then
$R^{(n+1)}\equiv R$. Hence $\omega^{(n)}<R$.
\end{proof}

\section{Incomparability between the jump operators}
In this section we show that the saturation jump operator and the halting jump operator
produce incomparable ceers in a very strong sense. A typical instance is still the halting
equivalence relation.

\begin{theorem}
For any $n\in\NN$, $H^+\not\leq H^{(n)}$.
\end{theorem}
\begin{proof}
For $n=0$ this is the content of Theorem 7.7. Assume $n\geq 1$ and assume toward
a contradiction $H^+\leq H^{(n)}$ via $f$. Let $s(x)$ be a one-one computable function
such that $\kappa(s(x))\!\downarrow\,=x$ for all $x\in\NN$. Let $h(x)$ be a one-one computable
function such that $\kappa(h(x))=\kappa(x)+n$ for all $x\in\NN$. For $x=\{ x_1,\dots,x_m\}$
define $g(x)=\{ s(0),\dots,s(n-1),h(x_1),\dots,h(x_m)\}$. Then $xH^+y\Leftrightarrow
g(x)H^+g(y)\Leftrightarrow f(g(x))H^{(n)}f(g(y))$ for any $x,y\in\NN^{<\omega}$.

We claim that for any $x\in\NN^{<\omega}$, $\kappa^n(f(g(x))\!\downarrow\,$. Granting the claim,
we would have $f(g(x))H^{(n)}f(g(y))\Leftrightarrow \kappa^n(f(g(x)))H\kappa^n(f(g(y))$,
which implies that $H^+\leq H$, a contradiction. To establish the claim, we prove a
stronger statement that for $k\leq n-1$ and any $x\in\NN^{<\omega}$, 
$\kappa^{k+1}(f(\{ s(0),\dots,s(k)\}\cup x))\!\downarrow\,$. 

The proof of the statement is by induction on $k\leq n-1$. For $k=0$ we need to
see that $\kappa(f(\{ s(0)\}\cup x))\!\downarrow\,$ for all $x$. Assume not, there would
be some $x_0$ such  that $[f(\{ s(0)\}\cup x_0)]_{H^{(n)}}$ is a singleton and therefore $[\{ s(0)\}\cup x_0]_H$ is computable, which is impossible. In general suppose the
statement is true for $k$ and $k+1\leq n-1$. Then for any $x$, $\kappa^{k+1}(f(\{s(0),\dots,
s(k+1)\}\cup x))\!\downarrow\,$ by the inductive hypothesis. Suppose for some $x_0$
$\kappa^{k+2}(f(\{ s(0),\dots,s(k+1)\}\cup x_0))\!\uparrow\,$. Without loss of generality
assume $[x_0]_H\cap [\{ s(0),\dots,s(k+1)\}]_H=\emptyset$. Then for any $y\in\NN$,
$$\begin{array}{rcl}
\kappa(y)\!\downarrow\,=k+1 &\Leftrightarrow&
(\{ s(0),\dots,s(k+1)\}\cup x_0)\, H^+\,(\{ s(0),\dots,s(k),y\}\cup x_0) \\
&\Leftrightarrow & f(\{ s(0),\dots,s(k+1)\}\cup x_0)\, H^{(n)}\,f(\{ s(0),\dots,s(k),y\}\cup x_0)
\\
& \Leftrightarrow & \kappa^{k+1}(f(\{ s(0),\dots,s(k+1)\}\cup x_0))= \\
& & \kappa^{k+1}(f(\{ s(0),\dots,s(k),y\}\cup x_0))
\end{array}
$$
Since the computations on the last line are all convergent, this gives a decision procedure
for $K_{k+1}=\{ y\,|\, \kappa(y)\!\downarrow\,=k+1\}$, a contradiction.
\end{proof}

On the other hand, it is also true that $H'\not\leq H^{n+}$ for any $n\in\NN$.
In fact we demonstrate a slightly stronger result.

\begin{theorem}
$H'\not\leq H^{\omega+}$.
\end{theorem}
\begin{proof}
Instead of considering the formal definition of $H^{\omega+}$ it is more convenient
to view the universe of $H^{\omega+}$ as the set of all finite trees with terminal nodes
labeled by natural numbers. For such a tree $t$ denote its set of labels for terminal nodes by $N_t$.

Now the proof uses a similar argument as that of Theorem 7.8. 
Assume $H'\leq H^{\omega+}$ via $f$. Let $s(x)$ be a computable
function such that $\kappa(s(x))\!\downarrow\,= x$, for any $x\in\NN$. Then for any $x,y\in\NN$,
$xHy\Leftrightarrow s(x)H's(y)\Leftrightarrow f(s(x))H^{\omega+}f(s(y))$.
Now for any $x\in\NN$, let $A_x=N_{f(s(x))}\cap K$ and $B_x=N_{f(s(x))}\setminus A_x$.
We have that $B_x$ is a constant finite set for all $x\in\NN$, because otherwise
it would follow that for some $x_0 \neq y_0$ the sets $[s(x_0)]_{H'}$ and $[s(y_0)]_{H'}$
are computably separable, which is impossible. It then follows that it is computable
to decide whether $f(s(x))H^{\omega+}f(s(y))$, which in turn implies that $H$ is 
computable, a contradiction.
\end{proof}

\begin{corollary}
$H^{\omega+}$ is not universal for all ceers.
\end{corollary}

There is an interesting contrast between Corollary 9.3 and Theorem 8.13.
\section{Open problems}
To motivate further research let us mention the following problems which seem to
be unsolved.

\begin{problem}
Is being universal a $\Sigma_3$-complete property among ceers?
\end{problem}

\begin{problem}
Is there a non-universal ceer $E$ with $E\equiv E'$?
\end{problem}

There are other, more techinical problems we encountered in our research. 
But the above ones seem to be the most significant problems we know of.
Also, in this paper we have been focusing on the general theory of reducibility
among ceers, not paying too much attention to the applications of this theory.
Most of the natural examples of ceers arise in algebra and logic. For example,
the word problems on finitely presented algebraic structures are essentially
ceers. The classical unsolvability results usually deal with a single equivalence
class. New problems emerge when they are investigated as ceers and when the strong
reducibility notion is involved. As an example let us mention the following problem.

\begin{problem}
Is the following statement true: for any ceer $R$ there is a finitely presented
semigroup whose word problem as a ceer is bi-reducible to $R$?
\end{problem}

It is obvious that a similar statement for finitely presented groups is false, since
if the word problem for a group is CC then it has to be computable.

Finally let us mention that the idea of considering reducibility between
equivalence relations as a framework for the study of complexity of classification
problems can also be applied elsewhere. For example, we could consider feasibly
computable reductions between computable equivalence relations and develop
a parallel theory as did in this paper. Such a theory would probably shed light
on structural complexity theory and provide deeper understanding of the open
questions in the field.

\end{document}